\documentclass[12pt]{article}
\usepackage{amsmath}
\usepackage{amssymb}

\usepackage[cp1251]{inputenc}
\usepackage[russian]{babel}
\begin{document}

\vskip 12pt\noindent
MSC 34C15, 34L30, 37C99

\vskip 12pt
\centerline {\bf Global solvability criteria for}
 \centerline {\bf some classes of nonlinear second order ordinary}
\centerline {\bf differential equations}
\vskip 12pt
\centerline {\bf G. A. Grigorian}
\centerline{ 0019,  Armenia,  c. Erevan, str. Bagramian, 24/5,}
\centerline{Institute of mathematics NAS of Armenia}

\centerline{E- mail: mathphys2@instmath.sci.am}

\vskip 12pt
Abstract. The  Riccati equation method is used to establish some global  solvability criteria  for some classes of second order  nonlinear ordinary differential equations. Two oscillation theorems are proved. The results are applied to the Emden - Fowler equation and  to the  Van der Pol  type equation.

\vskip 12pt

\vskip 12pt
{\bf \centerline { \S 1. Introduction}}
\vskip 12pt

Due to numerous applications the nonlinear ordinary differential equations occupy an important place in the theory of differential equations and  numerous works are devoted to them ([1 - 8] and cited works therein). Except in rare cases these equations cannot be integrated explicitly. Therefore, an important role plays   study of conditions of their global solvability, qualitative study of their solutions (e.g.  such characteristics of their solutions  as oscillation,  asymptotic behavior, stability and so on).

Let $p_0(t;w), \phantom{a} q_0(t;w)$ and $r_0(t;w)$  be continuous on  $[t_0;+\infty)\times(-\infty;+\infty)$ real valued functions, and let  $p_0(t;w) > 0,\phantom{a} t\ge t_0,\phantom{a} w\in(-\infty;+\infty)$. Consider the equation $$
(p_0(t;\phi(t))\phi'(t))' + q_0(t;\phi(t))\phi'(t) + r_0(t;\phi(t))\phi(t) = 0,\phantom{aaa} t\ge t_0. \eqno (1.1)
$$

Like the linear differential equations of the second order, this equation can be interpreted as a system of nonlinear differential equations of the first order (see [3], p. 381):
$$
\left\{
\begin{array}{l}
\phi'(t)= \frac{\psi(t)}{p_0(t;\phi(t))};\\
\phantom{a}\\
\psi'(t)= - r_0(t;\phi(t))\phi(t) - \frac{q_0(t;\phi(t))}{p_0(t;\phi(t))}\psi(t).
\end{array}
\right.
\eqno (1.2)
$$
By Peano's theorem (see [3], p. 21, 22) for every
$\phi_{(0)}$  and  $\phi_{(1)}$  (here and henceforth   $\phi_{(0)}$  and $\phi_{(1)}$  are real numbers)  and  $t_1\ge t_0$  the system  (1.2) has a solution  $(\phi(t), \psi(t))$ in the neighborhood of the point  $t_1$  (in the case  $t_1=t_0$  in some right  neighborhood  of  $t_0$),
satisfying the initial conditions: $\phi(t_1) = \phi_{(0)},\phantom{a}  \psi(t_1) = \phi_{(1)}$.
Therefore, for any   $\phi_{(0)},\phantom{a} \phi_{(1)}$   and   $t_1\ge t_0$  eq. (1.1) has a solution
in some neighborhood of the point $t_1$, satisfying the initial conditions:
$\phi(t_1) = \phi_{(0)}, \phantom{a} \phi'(t_1) = \phi_{(1)}$.

{\bf Remark 1.1}. The solution  $\phi(t)$ of eq.  (1.1), satisfying the initial conditions:   $\phi(t_1) = \\ = \phi_{(0)}, \phantom{a} \phi'(t_1) = \phi_{(1)}$ in general, is not the unique. However, under additional restrictions on the functions   $p_0(t;w), q_0(t;w)$ and $r_0(t;w)$ it is the  unique.
For example, due to (1.2),  if the functions  $f_1(t;u;v)\equiv \frac{v}{p_0(t;u)}, \phantom{a}f_2(t;u;v)\equiv r_0(t;u)u + \frac{q_0(t;u)}{p_0(t;u)} v$  satisfy the Lipschitz condition jointly $u, v$
in the region  $O\equiv \{(t;u;v): |t - t_1| \le \delta, \phantom{a}|u - \phi_{(0)}| \le M, \phantom{a}|v  - \phi_{(1)}| \le\\ \le N\}, \phantom{a}\delta >0,\phantom{a} M > 0,\phantom{a} N >0$ (see [3]. p. 13), then by virtue of the Picard - Lindellef's theorem (see [3], p. 19) the solution $\phi(t)$ of eq. (1.1) with
$\phi(t_1) = \phi_{(0)},\phantom{a} \phi'(t_1) = \phi_{(1)}$ exists on the interval $[t_1;t_2]$ and is the
unique, where $t_2\equiv \min\{\delta, \frac{\sqrt{M^2 + N^2}}{M_0}\},
 M_0\equiv \\ \equiv \max\limits_{(t;u;v) \in O} \sqrt{f_1^2(t;u;v) + f_2^2(t;u;v)}$
(as far as the solution $(\phi(t), \psi(t))$
of the system (1.2) with $\phi(t_1) = \phi_{(0)},\phantom{a} \psi(t_1) = \phi_{(1)}$ exists on the interval $[t_1;t_2]$).

{\bf Example 1.1}. For  $p_0(t;u) \equiv |u|^\sigma, \sigma > 0,\phantom{a} q_0(t;u)\equiv 0,\phantom{a} r_0(t;u) \equiv |u|^\nu,\phantom{a} \nu < 0$ the functions $f_j(t;u;v),\phantom{a} j=1,2,$ satisfy the Lipschtz condition in the region $|t - t_1| \le \delta,\phantom{a} \delta >\\> 0, \phantom{a}|u - 1| \le \frac{1}{2},\phantom{a} |v| \le 1$, but not in the region $|t - t_1| \le \delta, \phantom{a}\delta > 0,\phantom{a} |u - 1| \le 2, \phantom{a}|v| \le 1$.

{\bf Example 1.2}. For $p_0(t;u) \equiv 1 + t^2 + u^4, \phantom{a}q_0(t;u) \equiv t + u^3,\phantom{a} r_0(t;u) \equiv t^3 - u$ the functions
$f_j(t;u;v), \phantom{a}j=1,2,$ satisfy the Lipschtz condition in the region $|t - t_1| \le \delta, \phantom{a}|u| \le \\\le M, \phantom{a}|v| \le N,  \phantom{a}\delta > 0, \phantom{a}M > 0, \phantom{a}N > 0$.

In this paper the Riccati equation method is applied to establish some global existence criteria for Eq. (1.1). Two oscillatory theorems are proved. The obtained  results are applied to Emden - Fowler's equation and to the Van der Pol type equation.

\vskip 12pt
{\bf \centerline { \S 2.Auxiliary propositions}}
\vskip 12pt

Throughout of this paragraph we will assume that $\phi_0(t)$ is a solution of Eq.(1.1) on the interval $[t_0;T) \phantom{a}(T\le +\infty)$; $y_0(t) \equiv p_0(t;\phi_0(t)) \frac{\phi'_0(t)}{\phi_0(t)}$. Consider the Riccati equation
$$
y'(t) + \frac{y^2(t)}{p_0(t;\phi_0(t))} + \frac{q_0(t;\phi_0(t))}{p_0(t;\phi_0(t))} y(t) + r_0(t;\phi_0(t)) = 0, \phantom{aaa}t\in [t_0;T). \eqno (2.1)
$$
Assume $\phi_0(t)\ne 0, \phantom{a}t\in [t_1;t_2) (\subset [t_0;T))$.
It is not difficult to check, that the function  $y_0(t)$  is a solution of Eq. (2.1) on the interval  $[t_1;t_2)$. We have:
$$
\phi_0(t)=\phi_0(t_1)\exp\biggl\{\int\limits_{t_1}^t\frac{\phi'_0(\tau)}{\phi_0(\tau)}d\tau\biggr\}
=\phi_0(t_1)\exp\biggl\{\int\limits_{t_1}^t\frac{y_0(\tau)}{p_0(\tau;\phi_0(\tau))}d\tau\biggr\},\phantom{aaa} t\in [t_1;t_2). \eqno (2.2)
$$
Consider the linear equation
$$
z'(t) + \frac{y_0(t) + q_0(t;\phi_0(t))}{p_0(t;\phi_0(t))} z(t) + r_0(t;\phi_0(t)) = 0,\phantom{aaa} t\in [t_1;t_2). \eqno (2.3)
$$
By virtue of the Cauchy's formula the general solution of this equation  on the interval $[t_1;t_2)$  is given by the formula
$$
z(t)=c\exp\biggl\{- \int\limits_{t_1}^t \frac{y_0(\tau) + q_0(\tau;\phi_0(\tau))}{p_0(\tau;\phi_0(\tau))}d\tau\biggr\} -\phantom{aaaaaaaaaaaaaaaaaaaaaaaaaaaaaaaaaaaaaaa}
$$
$$
\phantom{aaaaaaaaa} - \int\limits_{t_1}^t\exp\biggl\{- \int\limits_{\tau}^t \frac{y_0(s) + q_0(s;\phi_0(s))}{p_0(s;\phi_0(s))}ds\biggr\}r_0(\tau;\phi_0(\tau)) d \tau,\phantom{aa} c=const. \eqno (2.4)
$$
By  (2.1) $y_0(t)$  is a solution of (2.3). Therefore from (2.4) it  follows:
$$
y_0(t)=y_0(t_1)\exp\biggl\{- \int\limits_{t_1}^t \frac{y_0(\tau) + q_0(\tau;\phi_0(\tau))}{p_0(\tau;\phi_0(\tau))}d\tau\biggr\} -\phantom{aaaaaaaaaaaaaaaaaaaaaaaaaaaaaaaaaaaaaaa}
$$
$$
\phantom{aaaaaaaaa} - \int\limits_{t_1}^t\exp\biggl\{- \int\limits_{\tau}^t \frac{y_0(s) + q_0(s;\phi_0(s))}{p_0(s;\phi_0(s))}ds\biggr\}r_0(\tau;\phi_0(\tau)) d \tau,\phantom{aaa} t\in [t_1;t_2). \eqno (2.5)
$$
By  (1.1) the equality
$$
[p_0(t;\phi_0(t))\phi'_0(t)]' + \frac{q_0(t;\phi_0(t))}{p_0(t;\phi_0(t))}[p_0(t;\phi_0(t))\phi_0'(t)] + r_0(t;\phi_0(t))\phi_0(t)=0,\phantom{aaa} t\ge t_0.
$$
is satisfied. Therefore
$$
p_0(t;\phi_0(t))\phi'_0(t) = p_0(t_1;\phi_0(t_1))\phi'_0(t_1)\exp\biggl\{- \int\limits_{t_1}^t \frac{q_0(\tau;\phi_0(\tau))}{p_0(\tau;\phi_0(\tau))}d\tau\biggr\} -\phantom{aaaaaaaaaaaaaaaaaaaaaaaaaaaaaaaaaaaaaaa}
$$
$$
\phantom{aaaaaaaaa} - \int\limits_{t_1}^t\exp\biggl\{- \int\limits_{\tau}^t \frac{q_0(s;\phi_0(s))}{p_0(s;\phi_0(s))}ds\biggr\}r_0(\tau;\phi_0(\tau))\phi_0(\tau) d \tau,\phantom{aaa} t\in [t_0;T). \eqno (2.6)
$$
Dividing both sides of this equality on  $p_0(t;\phi_0(t))$  and integrating from  $t_1$ to $t$ we obtain:
$$
\phi_0(t)=\phi_0(t_1) + p_0(t_1;\phi_0(t_1))\phi'_0(t_1)\int\limits_{t_1}^t\exp\biggl\{- \int\limits_{t_1}^{\tau} \frac{q_0(s;\phi_0(s))}{p_0(s;\phi_0(s))}d s \biggr\}\frac{d\tau}{p_0(\tau;\phi_0(\tau))}-
$$
$$
-\int\limits_{t_1}^t\frac{d\tau}{p_0(\tau;\phi_0(\tau))}\int\limits_{t_1}^\tau\exp\biggl\{- \int\limits_{s}^\tau \frac{q_0(\xi;\phi_0(\xi))}{p_0(\xi;\phi_0(\xi))}d\xi\biggr\}r_0(s;\phi_0(s))\phi_0(s) d s. \eqno (2.7)
$$

Let  $p_1(t;w),\phantom{a} q_1(t;w)$ and  $r_1(t;w)$ be real valued continuous functions on the $[t_0;+\infty)\times(-\infty;+\infty)$, and let  $p_1(t;w) > 0,\phantom{a} t\ge 0,\phantom{a} w\in(-\infty;+\infty)$. Along with the (1.1) consider the  equation
$$
(p_1(t;\phi(t))\phi'(t))' + q_1(t;\phi(t))\phi'(t) + r_1(t;\phi(t))\phi(t) = 0,\phantom{aaa} t\ge t_0. \eqno (2.8)
$$
Throughout of this paragraph we will assume that  $\phi_1(t)$ is a solution of Eq.(2.8) on the interval $[t_0;T) \phantom{a}(T\le  +\infty)$; $y_1(t) \equiv p_1(t;\phi_1(t)) \frac{\phi'_1(t)}{\phi_1(t)}$.
 Consider the Riccati equation
$$
y'(t) + \frac{y^2(t)}{p_1(t;\phi_1(t))} + \frac{q_1(t;\phi_1(t))}{p_1(t;\phi_1(t))} y(t) + r_1(t;\phi_1(t)) = 0,\phantom{aaa} t\in [t_0;T). \eqno (2.9)
$$
Let  $\phi_1(t)\ne0, \phantom{a}t\in [t_1;t_2)\phantom{a}(\subset [t_0;T)$. Then as in the case of Eq. (2.1) the function
$y_1(t)$
is a solution of Eq. (2.9) on the interval $[t_1;t_2)$  and
$$
\phi_1(t)
=\phi_1(t_1)\exp\biggl\{\int\limits_{t_1}^t\frac{y_1(\tau)}{p_1(\tau;\phi_1(\tau))}d\tau\biggr\},\phantom{aaa} t\in [t_1;t_2). \eqno (2.10)
$$
Since $y_0(t)$ and $y_1(t)$ are solutions of Eq. (2.1)  and   Eq. (2.9) respectively,  we have:
$$
y_j'(t) + \frac{y_j^2(t)}{p_j(t;\phi_j(t))} + \frac{q_j(t;\phi_j(t))}{p_j(t;\phi_1(t))} y_j(t) + r_j(t;\phi_j(t)) = 0, \phantom{aaa} j=0,1.
$$
Therefore,
$$
[y_1(t) - y_0(t)]' + \frac{y_0(t) + y_1(t) + q_{1-j}(t;\phi_{1-j}(t))}{p_{1-j}(t;\phi_{1-j}(t))}[y_1(t) - y_0(t)] +\phantom{aaaaaaaaaaaaaaaaaaaaaaaaaa}
$$
$$
+ \biggl[\frac{1}{p_1(t;\phi_1(t))} - \frac{1}{p_0(t;\phi_0(t))}\biggr]y_j^2(t) + \biggl[\frac{q_1(t;\phi_1(t))}{p_1(t;\phi_1(t))} - \frac{q_0(t;\phi_0(t))}{p_0(t;\phi_0(t))}\biggr] y_j(t) +\phantom{aaaaaaa}
$$
$$
\phantom{aaaaaaaaaaaaaaa}+ r_1(t;\phi_1(t)) - r_0(t;\phi_0(t)) = 0,\phantom{aaa} t\in [t_1;t_2), \phantom{a} j=0,1. \eqno (2.11)
$$
By  (2.5) from  (2.11) we have:
$$
y_1(t) - y_0(t) = [y_1(t_1) - y_0(t_1)]\exp\biggl\{- \int\limits_{t_1}^t\frac{y_0(\tau) + y_1(\tau) + q_{1-j}(\tau;\phi_{1-j}(\tau))}{p_{1-j}(\tau;\phi_{1-j}(\tau))}d\tau\biggr\} -
$$
$$
- \int\limits_{t_1}^t\exp\biggl\{-\int\limits_{\tau}^t\frac{y_0(s) + y_1(s) + q_{1-j}(s;\phi_{1-j}(s))}{p_{1-j}(s;\phi_1(s))}ds\biggr\}\biggl[\biggl(\frac{1}{p_1(\tau;\phi_1(\tau))} - \frac{1}{p_0(\tau;\phi_0(\tau))}\biggr) y_j^2(\tau) + \biggr.
$$
$$
\biggl. +\biggl(\frac{q_1(\tau;\phi_1(\tau))}{p_1(\tau;\phi_1(\tau))} - \frac{q_0(\tau;\phi_0(\tau))}{p_0(\tau;\phi_0(\tau))}\biggr) y_j(\tau) + r_1(\tau;\phi_1(\tau)) - r_0(\tau;\phi_0(\tau))\biggr]d\tau, \phantom{a}j=0,1. \eqno (2.12_j)
$$

{\bf Lemma 2.1}. {\it Assume $y_0(t_1) \ge 0$,  and  the inequality
$$
r_0(t;w) \le 0\phantom{aa} \mbox{for} \phantom{aa}  0 < w \le Y_0(t),\phantom{aa} t\in [t_1;t_2), \eqno (2.13)
$$
is satisfied  where  $Y_0(t)\equiv\max\limits_{\xi\in [t_1;t]}|\phi_0(\xi)|$. Then
$$
y_0(t)\ge 0, \phantom{aaa}  t\in [t_1;t_2).   \eqno (2.14)
$$
Moreover  if  $y_0(t_1) > 0$, then
 $$
y_0(t)> 0, \phantom{aaa}  t\in [t_1;t_2),   \eqno (2.15)
$$
}

Proof. It follows from (2.13), that
$r_0(t,\phi(t)) \le 0, \phantom{a} t\in [t_1;t_2)$. Then
$$
\int\limits_{t_1}^t\exp\biggl\{- \int\limits_\tau^t\frac{y_0(s) + q_0(s,\phi_0(s))}{p_0(s,\phi_0(s))}d s\biggr\}r_0(\tau,\phi_0(\tau))d \tau \le 0, \phantom{a} t\in [t_1;t_2).
$$
From here, from (2.5) and from the inequality $y_0(t_1) \ge 0$ ($y_0(t_1) > 0$) it follows (2.14) ((2.15)).
The lemma is proved.

{\bf Definition 2.1}. {\it  The set  $[t_1;t_2)$ is called the
 maximum  existence interval for the solution
 $\phi_0(t)$ $(y_0(t), \phantom{a}\phi_1(t),\phantom{a} y_1(t))$,
if  $\phi_0(t)$ $(y_0(t),\phantom{a} \phi_1(t),\phantom{a} y_1(t))$ exists on the interval $[t_1;t_2)$
and cannot be  continued to right from $t_2$ as a solution of Eq. (1.1)  ((2.1),   (2.8), (2.9)).}

{\bf Lemma 2.2}. {\it Assume   $y_0(t)\ge 0, \phantom{a}  t\in [t_1;t_2)$, and the inequality
$$
\int\limits_{t_1}^{t_2}\frac{y_0(\tau)}{p_0(\tau;\phi_0(\tau))}d\tau < +\infty. \eqno (2.16)
$$
is fulfilled. Then  $[t_1;t_2)$  is not the maximum  existence interval for $y_0(t)$.}

Proof. Since  $y_0(t)$ is nonnegative,   it follows from (2.2) and  (2.16) that there exists a finite limit
$$
\lim\limits_{t \to t_2 -0}\phi_0(t)=\phi_0(t_1)\exp\biggl\{\lim\limits_{t \to t_2 -0}\int\limits_{t_1}^t\frac{y_0(\tau)}{p_0(\tau;\phi_0(\tau))}d\tau\biggr\} \ne 0. \eqno (2.17)
$$
By  (2.6) from here  and from the continuity of  $p_0(t;w), \phantom{a}q_0(t;w),\phantom{a} r_0(t;w)$
it follows the existence of a finite limit  $\lim\limits_{t\to t_2 - 0}\phi'_0(t)$.
Therefore, $[t_1;t_2)$  is not the maximum  existence interval for  $\phi_0(t)$, and hence
by  (2.17) the function  $p_0(t;\phi_0(t))\frac{\phi'_0(t)}{\phi_0(t)}$ is defined on the interval $[t_1;t_2 + \varepsilon)$   for some  $\varepsilon > 0$. It follows from here that  $[t_1;t_2)$ is not the maximum existence interval for $y_0(t)$. The proof of the lemma is completed.

Let   $u=u(t) \phantom{a} (\ne0),\phantom{a} v=v(t),\phantom{a} x=x(t),$  $[t_0;+\infty), \phantom{a} P(t),  Q(t)$  and   $R(t)$  be real valued continuous functions on the  interval  $[t_0;+\infty)$, and let $P(t) > 0,\phantom{a} t \ge t_0$.   Assume
$$
I^+_{u,v}(t_1;t)\equiv \int\limits_{t_1}^t\exp\biggl\{- \int\limits_{t_1}^\tau v(s)d s \biggr\}\frac{d\tau}{u(\tau)},\phantom{aaa} I^-_{v,x}(t_1;t) \equiv \int\limits_{t_1}^t\exp\biggl\{- \int\limits_\tau ^t v(s)d s \biggr\}x(\tau)d\tau,
$$
Set:
$$
F(t_1;t;c_1;c_2) \equiv |c_1|\exp\biggl\{ c_2 I^+_{P,Q}(t_1;t) - \int\limits_{t_1}^t I^-_{Q,R}(t_1;\tau)\frac{d\tau}{P(\tau)}\biggr\},\phantom{a} t_1,\phantom{a} t \ge t_0.
$$

{\bf Lemma 2.3}. {\it Assume  $y_0(t)\ge0,\phantom{a} t\in [t_1;t_2)$, and  for some  $\varepsilon > 0$ the following inequalities are satisfied:
$$
p_0(t;w) \ge P(t), \phantom{a} \frac{q_0(t;w)}{p_0(t;w)} \ge Q(t), \phantom{a} R(t) \le r_0(t;w) \le 0   \eqno (2.18)
$$
for  $|w| \le F(t_1;t;c_1;c_2) + \varepsilon, \phantom{a}t\in [t_1;t_2)$, where $c_1\equiv \phi_0(t_1)\ne 0,\phantom{a} c_2 \equiv y_0(t_1)$.
Then
$$
\int \limits_{t_1}^t\frac{y_0(\tau)}{p_0(\tau;\phi_0(\tau))}d\tau\le y_0(t_1)I^+_{P, Q}(t_1;t) - \int\limits_{t_1}^t I^-_{Q, R}(t_1;\tau)\frac{d\tau}{P(\tau)},\phantom{aaa} t\in [t_1;t_2). \eqno (2.19)
$$
}

Proof. Suppose  for some  $t_3\in (t_1;t_2)$  the inequality (2.19) is false. Then because this inequality holds for $t=t_1$, then taking into account (2.2) we get that there exists $t_4\in~ (t_1;t_3)$  such, that
$$
|\phi_0(t_4)| > F(t_1;t_4; c_1; c_2); \eqno (2.20);
$$
$$
|\phi_0(t)| \le F(t_1;t; c_1; c_2) + \varepsilon, \phantom{aaa}t\in [t_1;t_4).
$$
From the last inequality and from (2.20) it follows, that
$$
p_0(t;\phi_0(t)) \ge P(t),\phantom{a} \frac{q_0(t;\phi_0(t))}{p_0(t;\phi_0(t))}\ge Q(t); \phantom{a}R(t) \le r_0(t;\phi_0(t)) \le 0,\phantom{aaa} t\in [t_1;t_4].
$$
By virtue of nonnegativity of  $y_0(t)$ from here and  from (2.5) it follows:
$$
\frac{y_0(t)}{p_0(t;\phi_0(t))} \le \frac{c_2}{P(t)}\exp\biggl\{ - \int\limits_{t_1}^t Q(\tau)d\tau\biggr\}
 - \frac{1}{P(t)} I^-_{Q,R}(t_1;t);\phantom{aaa} t\in [t_1;t_4].
$$
Integrating this inequality from  $t_1$  to  $t_4$ we obtain:
$$
\int\limits_{t_1}^{t_4}\frac{y_0(\tau)}{p_0(\tau;\phi_0(\tau))} d \tau \le c_2 I^+_{P, Q}(t_1;t_4) - \int\limits_{t_1}^{t_4}I^-_{Q, R}(t_1;\tau)\frac{d\tau}{P(\tau)}.
$$
Consequently, $|\phi_0(t_4)| \le F(t_1;t_4;c_1;c_2)$ which contradicts (2.20). The obtained contradic-\\tion proves (2.19). The lemma is proved.

Set $G_{x}(t_1;t;c_1;c_2)\equiv |c_1|\exp\biggl\{c_2 I^+_{P, Q}(t_1;t) + \int\limits_{t_1}^t\frac{x(\tau)}{P(\tau)}d\tau\biggr\}$.

{\bf Lemma 2.4}. {\it Let  $y_0(t) \ge 0,\phantom{a} t\in [t_1;t_2)$,
and let for some  $\varepsilon >0$  the following inequalities be satisfied:
$$
p_0(t;w) \ge P(t), \phantom{a}q_0(t;w) \ge Q(t)\ge 0,\phantom{a} r_0(t;w) \le 0,\phantom{a} \left|p_0(t;w)\frac{r_0(t;w)}{q_0(t;w)}\right|\le \widetilde{Q}(t) \phantom{aa}\mbox{for}\phantom{aa} |w|\le
$$
$\le G_{M}(t_1;t;c_1;c_2) + \varepsilon ,\phantom{a} t\in [t_1;t_2)$, where  $\widetilde{Q}(t)$  is a continuous function on the interval $[t_1;t_2)$, $M(t)\equiv \max\limits_{\xi\in[t_1;t]}\{\widetilde{Q}(\xi)\}, \phantom{a}c_1 = \phi_0(t_1),\phantom{a} c_2 = y_0(t_1)$.   Then
$$
\int\limits_{t_1}^t\frac{y_0(\tau)}{p_0(\tau;\phi_0(\tau))} d \tau \le y_0(t_1) I^+_{P, Q}(t_1;t) - \int\limits_{t_1}^t\frac{M(\tau)}{P(\tau)} d\tau,\phantom{aaa} t\in [t_1;t_2).
$$
}
This lemma can be proved by analogy of Lemma 2.3. In its proof can be used the following  easily verifiable inequality
$$
\int\limits_{t_1}^t\exp\biggl\{- \int\limits_\tau ^t \frac{y_0(s) + q_0(s;\phi_0(s))}{p_0(s;\phi_0(s))}\biggr\}|r_0(\tau;\phi_0(\tau))|d\tau \le M(t),\phantom{aaa} t\in [t_1;t_2).
$$

{\bf Lemma 2.5}. {\it Let the following conditions be satisfied:

\noindent
a$_1) \phantom{a} y_1(t_1) > y_0(t_1)$;

\noindent
b$_1)  \phantom{a} \phi_1(t_1) \ge \phi_0(t_1) > 0, \phantom{aaa}  \phi'_1(t_1) > \phi'_0(t_1)\ge 0$  or

$\phi_1(t_1) \le \phi_0(t_1) < 0, \phantom{a}  \phi'_1(t_1) < \phi'_0(t_1)\le 0$;

\noindent
c$_1)\phantom{a} p_0(t;w) \le p_1(t;w_1), \phantom{a}r_0(t;w) \ge r_1(t;w_1)$, for  $t\in [t_1;t_2),\phantom{a} |w|\le |w_1|,\phantom{a} w,\phantom{a}$\\
\phantom{aaaa}$w_1 \in(-\infty; +\infty)$;

\noindent
d$_1)\phantom{a}  p_0(t;w) \ge p_1(t;w_1)$  for   $t\in [t_1;t_2), |w|\ge |w_1|, w, w_1 \in(-\infty; +\infty)$;

\noindent
e$_1)\phantom{a} \frac{q_0(t;w)}{p_0(t;w)}\le  \frac{q_1(t;w_1)}{p_1(t;w_1)}$; $y_1(t) \ge0$  or   $y_0(t) \ge 0$   for  $t\in [t_1;t_2), |w|\le |w_1|$,

 $w, w_1 \in(-\infty; +\infty)$.

\noindent
Then
$$
y_1(t) > y_0(t), \phantom{aaa}t\in [t_1;t_2). \eqno (2.21)
$$
}

Proof. Suppose (2.23) is false. Then from  a$_1)$ it follows:
$$
y_1(t) > y_0(t), \phantom{aaa}t\in [t_1;t_3); \eqno (2.22)
$$
$$
y_1(t_3)= y_0(t_3)\eqno (2.23)
$$
for some $t_3\in (t_1;t_2)$.  From   b$_1)$ it follows:
$$
|\phi_1(t)| \ge |\phi_0(t)|,\phantom{aaa} t\in [t_1;t_4),  \eqno (2.24)
$$
for some  $t_4\in[t_1;t_3]$. Let us show, that
$$
|\phi_0(t)| \le |\phi_1(t)|,\phantom{aaa} t\in [t_1;t_3). \eqno (2.25)
$$
Suppose,  it is not true. Then it follows from  (2.24), that
$$
|\phi_0(t)| \le |\phi_1(t)|,\phantom{aaa} t\in [t_1;t_5], \eqno (2.26)
$$
$$
|\phi_0(t)| > |\phi_1(t)|, \phantom{aaa}t\in (t_5;t_3),  \eqno (2.27)
$$
where $t_5 = \sup\{t \in [t_1;t_3): |\phi_0(t)| \le |\phi_1(t)|\} \in (t_1; t_3)$.  On the strength of   (2.2) and (2.10) from    d$_1)$, (2.23),  (2.26) and (2.27) it follows:
$$
|\phi_0(t)|=|\phi_0(t_5)|\exp\biggl\{\int\limits_{t_1}^t\frac{y_0(\tau)}{p_0(\tau;\phi_0(\tau))}d\tau\biggr\} \le |\phi_1(t_5)|\exp\biggl\{\int\limits_{t_1}^t\frac{y_1(\tau)}{p_1(\tau;\phi_1(\tau))}d\tau\biggr\}=|\phi_1(t)|,
$$
$t\in [t_5;t_3)$, which contradicts (2.27). The obtained contradiction proves (2.25).
By (2.12) from  c$_1)$,  e$_1)$ and  (2.25) it follows, that $y_1(t_3) > y_0(t_3)$,
which contradicts  (2.23). The obtained contradiction proves (2.21). The  lemma is proved.

{\bf Remark 2.1}. {\it It follows from the conditions c$_1)$ and d$_1)$, that $p_0(t;w) = p_1(t;w), \phantom{a} t\in\\ \in [t_1;t_2),\phantom{a} w\in (-\infty;+\infty)$ and $p_0(t;w)$ increases (in the wide sense) by $w$ on the interval $[0;+\infty)$ and decreases (in the wide sense) by $w$ on the interval $(-\infty;0]$ for every $t\in [t_1;t_2)$, in particular when $p_0(t;w) = p_1(t;w) = p(t),\phantom{a} t\in [t_1;t_2)$, then the conditions в$_1)$ and г$_1)$ are satisfied.}

Consider the Riccati equation
$$
y'(t) + \frac{y^2(t)}{P(t)} + \frac{Q(t)}{P(t)}y(t) + R(t) = 0,\phantom{aaa} t\ge t_0. \eqno (2.28)
$$

{\bf Lemma 2.6}. {\it Let the following conditions be satisfied:

\noindent
a$_2)\phantom{a} y_0(t) \ge 0, \phantom{aaa}t\in [t_1;t_2) (\subset [t_0;T))$;

\noindent
b$_2)\phantom{a} P(t) \ge p_0(t;\phi_0(t)), \phantom{a}\frac{Q(t)}{P(t)}\le \frac{q_0(t;\phi_0(t))}{p_0(t;\phi_0(t))}, \phantom{a}R(t) \le r_0(t;\phi_0(t)), \phantom{a}t\in [t_1;t_2)$.

\noindent
Then the solution  $y_2(t)$  of Eq.  (2.28), satisfying the condition  $y_2(t_1) \ge y_0(t_1)$,
 exists on the interval $[t_1;t_2)$.}

Proof. Let  $[t_1;t_3)$ be  the  maximum existence  interval for $y_2(t)$.
We must show  that  $t_3\ge t_2$. Suppose  $t_3 < t_2$. Then taking into account the inequality   $y_2(t_1) \ge y_0(t_1)$,
 and $(2.12_0)$ we conclude: it follows from  a$_2)$ and   b$_2)$, that   for   $p_1(t;w)\equiv P(t)$,  $q_1(t;w)\equiv Q(t)$,\phantom{a} $r_1(t;w)\equiv R(t)$ the inequality  $y_2(t) \ge y_0(t), \phantom{a}t\in [t_1;t_3)$,  holds.
Since   $y_0(t)$ is continuous on the interval $[t_1;t_3]$,  it follows from the last inequality,
that the function   $f(t)\equiv \int\limits_{t_1}^t\frac{y_2(\tau)}{P(\tau)}d\tau$
is bounded from below on the interval $[t_1;t_3)$. It follows from here (see [9, p. 3,  Lemma 2.2])  that $[t_1;t_3)$   is not the maximum existence interval for  $y_2(t)$. The obtained contradiction shows, that $t_3\ge t_2$.
The  lemma is proved.

{\bf Lemma 2.7}. {\it Assume  $r_0(t;w) \ge 0,\phantom{a} t\in [T_0;T),\phantom{a} w\in (-\infty;+\infty)  (t_0 \le T_0 < T \le~  +\infty), \linebreak \phi_0(T_1) = 0$ for some  $T_1\in [T_0;T)$, and   $\sup\limits_{t\in [T_0;T)} |\phi_0(t)| > 0$.  Then  $\phi_0(t)$ changes  sign on the interval $[T_0;T)$.}

Proof. Since  $\sup\limits_{t\in [T_0;T)} |\phi_0(t)| > 0$,
 there exists  $T_2\in [T_0;T)$ such, that $\phi_0(T_2)\ne 0$.  Let  $T_2< T_1$
(the proof in the case  $T_2>T_1$  by  analogy). Then there exists  $T_3\in (T_2;T_1]$ such, that
$$
\phi_0(T_3) = 0,  \eqno (2.29)
$$
$$
\phi_0(t)\ne 0,\phantom{aaa} t\in [T_2;T_3). \eqno (2.30)
$$
By Lagrange's mean value  theorem it follows from (2.29) the existence of a $\xi \in (T_2;T_3)$    such, that  $\phi_0'(\xi) = - \frac{\phi_0(T_2)}{T_3 - T_2}$. Hence
$$
sign \hskip 2pt \phi'_0(\xi) = - sign \hskip 2pt \phi_0(T_2).  \eqno (2.31)
$$
To complete the proof of the  lemma it is enough to show, that
$$
sign \hskip 2pt \phi'_0(T_3) = - sign \hskip 2pt \phi_0(T_2).  \eqno (2.32)
$$
It follows from (2.30) that  $y_0(t)$  exists on the interval  $[T_2;T_3)$. Therefore, by  (2.5) it follows from the nonnegativity of $r_0(t;w)$, that $sign
  \hskip 3pt y_0(T_3) = sign \hskip 3pt  y_0(\xi)$. From here and (2.31) it follows (2.32). The lemma is proved.

\vskip 12pt
{\bf \centerline { \S 3. Some  global solvability and oscillatory criteria}}
\vskip 12pt

Let $P(t),\phantom{a} Q(t)$ and $R(t)$ be the same functions as in the previous paragraph.

{\bf Theorem 3.1}. {\it Assume  $\phi_{(0)}\ne 0,\phantom{a} \frac{\phi_{(1)}}{\phi_{(0)}}\ge0$,
and let for some  $\varepsilon > 0$ the following  inequalities are satisfied
$$
p_0(t;w) \ge P(t),\phantom{a} \frac{q_0(t;w)}{p_0(t;w)} \ge Q(t),\phantom{a} R(t) \le r_0(t;w) \le 0 \eqno (3.1)
$$
for $|w| \le F\bigl(t_0; t; \phi_{(0)}; p(t_0;\phi_{(0)})\frac{\phi_{(1)}} {\phi_{(0)}}\bigr) + \varepsilon, \phantom{a}t\ge t_0$.   Then the solution  $\phi_0(t)$  of Eq. (1.1),  satisfying the initial value conditions: $\phi_0(t_0) = \phi_{(0)},\phantom{a} \phi_0(t_0) = \phi_{(1)}$, exists on the interval $[t_0;+\infty)$.
The function  $|\phi_0(t)|$ is positive, nondecreasing and satisfies the estimate
$$
|\phi_0(t)|\le F\biggl(t_0; t; \phi_{(0)}; p_0(t_0;\phi_{(0)})\frac{\phi_{(1)}} {\phi_{(0)}}\biggr),\phantom{aaa} t\ge t_0. \eqno (3.2)
$$
In this case if   $\phi_{(1)} \ne 0$, then
$$
\phi'(t) \ne 0,\phantom{aaa} t\ge t_0. \eqno (3.3)
$$
}

Proof. Suppose   $\phi_{(0)} >0$   (the proof in the case
$\phi_{(0)} <0$ by analogy), $\phi_0(t)$ is  the solution of Eq. (1.1),
satisfying the initial value conditions: $\phi_0(t_0) = \phi_{(0)},\phantom{a} \phi'_0(t_0) = \phi_{(1)}$ (existence of $\phi_0(t)$ follows from the connection between (1.1) and (1.2) and from the Peano's theorem; see [3, p. 21, Theorem 2.1]).
Let  $[t_0;T)$  be  the maximum  existence interval for  $y_0(t)\equiv p_0(t;\phi_0(t))\frac{\phi'_0(t)}{\phi_0(t)}$  (it assumes, that  $\phi_0(t)$
exists on the interval $[t_0;T)$  and does not vanish on it). Show  that
$$
y_0(t)\ge 0,\phantom{aaa} t\in [t_0;T). \eqno (3.4)
$$
Suppose, that it is false.  Then since  $y_0(t_0)= p_0(t_0;\phi_{(0)})\frac{\phi_{(1)}}{\phi_{(0)}} \ge 0$ there exists $\widetilde{t}_0,\phantom{a} \widetilde{t}_1$  such  that
$$
y_0(t) \ge 0,\phantom{aaa} t\in [t_0;\widetilde{t}_0], \eqno (3.5)
$$
$$
y_0(t) < 0, \phantom{aaa}t\in (\widetilde{t}_0;\widetilde{t}_1). \eqno (3.6)
$$
By  (3.1)  for  $c_1= \phi_{(0)}, \phantom{a}c_2 = p_0(t_0;\phi_{(0)})\frac{\phi_{(1)}}{\phi_{(0)}},\phantom{a} t_1=t_0$ the inequalities
(2.20) of  Lemma~2.3 hold. Therefore taking into account (2.2) we  have
$$
\phi_0(t)\le F\biggl(t_0; t; \phi_{(0)}; p_0(t_0;\phi_{(0)})\frac{\phi_{(1)}} {\phi_{(0)}}\biggr),\phantom{aaa} t\in [t_0;\widetilde{t}_0].
$$
By  (2.2) from here and from (3.6) it follows
$$
\widetilde{M}(t)\equiv \max\limits_{\xi\in [t_0;t]}\phi_0(\xi)\le F\biggl(t_0; t; \phi_{(0)}; p_0(t_0;\phi_{(0)})\frac{\phi_{(1)}} {\phi_{(0)}}\biggr),\phantom{aaa} t\in [t_0;\widetilde{t}_1).
$$
Then taking into account the third  of the inequalities (3.1) we get  $r_0(t;w) \le 0$ for  $0< w \le  \widetilde{M}(t),\phantom{a} t\in  [t_0; \widetilde{t}_1)$.
On the strength of  Lemma 2.1 we conclude from here, that $y_0(t)\ge 0$ for  $t \in[t_0; \widetilde{t}_1)$ which  contradicts (3.6). The obtained contradiction proves (3.4).
Show, that $T= +\infty$.   Suppose  $T<+\infty$. By virtue of  Lemma 2.3  from  (3.1) and (3.3) it follows:
$$
\int\limits_{t_0}^t\frac{y_0(\tau)}{p_0(\tau;\phi_0(\tau))}d\tau \le y_0(t_0)I^+_{P,Q}(t_0;t) - \int\limits_{t_0}^t I^-_{Q, R}(t_0;\tau)\frac{d\tau}{P(\tau)},\phantom{aaa} t\in [t_0;T). \eqno (3.7)
$$
Therefore, $\int\limits_{t_0}^T\frac{y_0(\tau)}{p_0(\tau;\phi_0(\tau))}d\tau < +\infty$. On the strength of  Lemma 2.2 it follows from here that  $[t_0;T)$ is not the maximal existence interval for $y_0(t)$. The obtained contradiction shows that  $T=+\infty$. Hence (3.7) is valid for all $t\ge t_0$.  Therefore, by  (2.2)  the inequality (3.2) is valid.
Since $\phi_0(t_0) = \phi_{(0)} > 0$, and $y_0(t)$  is nonnegative  by  (2.2) the function  $\phi_0(t)$  is positive and nondecreasing on the interval $[t_0;+\infty)$.  And if  $\phi_{(1)} > 0$,
then by (2.2)  and  Lemma 2.1 the inequality (3.3) is fulfilled. The  theorem is proved.

{\bf Remark 3.1}. {\it  A solution $\phi_*(t)$  of the equation
$$
(P(t)\phi'(t))' + Q(t)\phi'(t) + R(t)\phi(t) = 0, \phantom{aaa}, t\ge t_o,
$$
such  that $\phi_*(t) \ne 0, \phantom{a} t\in [t_1;t_2)$, is connected with the function $F$ by the following relation
$$
|\phi_*(t_1)|\exp\bigg\{\frac{\phi_*(t)}{\phi_*(t_1)}\biggr\} = F\biggl(t_1;t;\phi_*(t_1);P(t_1)\frac{\phi_*'(t_1)}{\phi_*(t_1)}\biggr)
\exp\biggl\{1 -
 P(t_1)\frac{\phi_*'(t_1)}{\phi_*(t_1)}\times  \phantom{aaaaaaaaaaaaaaaaaaaaaaaaaaaaaaaaaaaa}
 $$

$$
\phantom{aaaaaaa} \times
 \int\limits_{t_1}^t\frac{d\tau}{P(\tau)}\int\limits_{t_1}^\tau \exp\biggl\{- \int\limits_s^\tau\frac{Q(\xi)}{P(\xi)} d \xi\biggr\} R(s) I^+_{P,Q}(t_1;s) ds + \frac{1}{\phi_*(t_1)}(K^2\phi_*)(t)\biggr\}, \phantom{aaa} \eqno (3.8)
$$
where $K$ - is the integral operator
$$
(K\phi_*)(t)\equiv \int\limits_{t_1}^t\frac{d\tau}{P(\tau)}\int\limits_{t_1}^\tau \exp\biggl\{- \int\limits_s^\tau\frac{Q(\xi)}{P(\xi)} d \xi\biggr\} R(s)\phi_*(s) d s, \phantom{aaa}, t\in [t_1;t_2).
$$
Indeed, since $y_*(t)\equiv P(t)\frac{\phi_*'(t)}{\phi_*(t)}$ is a solution of eq. (2.8) on the interval $[t_1;t_2)$,  by the Cauchy's formula
$$
y_*(t) = y_*(t_1)\exp\biggl\{ - \int\limits_{t_1}^t\frac{Q(\tau)}{P(\tau)} d \tau - \int\limits_{t_1}^t\frac{y_*(\tau)}{P(\tau)} d \tau\biggr\} -\phantom{aaaaaaaaaaaaaaaaaaaaaaaaaaaaaaaaaaaaa}
$$
$$
\phantom{aaaaaaaaaaaaaaaaaaaaaaaa} - \int\limits_{t_1}^t\exp\biggl\{- \int\limits_\tau^t\frac{Q(s)}{P(s)} d s - \int\limits_\tau^t\frac{y_*(s)}{P(s)} d s\biggr\} R(\tau) d \tau, \phantom{aaa} t\in [t_1;t_2)
$$
Multiplying both sides of this equality on $\frac{1}{P(t)}\exp\biggl\{\int\limits_{t_1}^t\frac{y_*(\tau)}{P(\tau)} d\tau\biggr\}$ and integrating from $t_1$ to $t$ taking into account the equality
$\phi_*(t) = \phi_*(t_1)\exp\biggl\{\int\limits_{t_1}^t \frac{y_*(\tau)}{P(\tau)} d\tau\biggr\}, \phantom{a} t\in[t_1;t_2)$, we obtain:
$$
\frac{\phi_*(t)}{\phi_*(t_1)} = 1 + P(t_1)\frac{\phi_*'(t_1)}{\phi_*(t_1)} I^+_{P,Q}(t_1;t) - \frac{1}{\phi_*(t_1)}(K\phi_*)(t), \phantom{aaa} t\in [t_1;t_2).
$$
After making the first iteration in this equality, taking its exponential  and multiplying by $|\phi_*(t_1)|$ we come to (3.8). The question of an applicability of the equality (3.8) for establishing effective criteria of global solvability of Eq (1.1) is an issue of separate study.}

Using  Lemma 2.4 in place of   Lemma 2.3 by analogy it can be proved

{\bf Theorem 3.2}. {\it Let  $\phi_{(0)}\ne 0, \phantom{a}\frac{\phi_{(1)}}{\phi_{(0)}}\ge0$,  and let for some $\varepsilon > 0$  the following inequalities be satisfied:
$$
p_0(t;w) \ge P(t),\phantom{a} q_0(t;w)\ge Q(t) \ge0,\phantom{a} r_0(t;w)\le 0, \phantom{a} \left|p_0(t;w)\frac{r_0(t;w)}{q_0(t;w)}\right|\le \widetilde{Q}(t)
$$
 for $|w|\le G_{M}\biggl(t_0;t;\phi_{(0)};p_0(t_0;\phi_{(0)})\frac{\phi_{(1)}}{\phi_{(0)}}\biggr) + \varepsilon, \phantom{a}t\ge t_0$. Then the solution $\phi_0(t)$  of Eq. (1.1),
satisfying the initial value conditions  $\phi_0(t_0) = \phi_{(0)},\phantom{a} \phi'(t_0)=\phi_{(1)}$,
exists on the interval $[t_0;+\infty)$. The function  $|\phi_0(t)|$ is positive, nondecreasing and satisfies the inequality
$$
|\phi_0(t)|\le G_{M}\biggl(t_0;t;\phi_{(0)};p_0(t_0;\phi_{(0)})\frac{\phi_{(1)}}{\phi_{(0)}}\biggr),\phantom{a} t\ge t_0.
$$
And if  $\phi_{(1)}\ne 0$, then  $\phi'_0(t)\ne 0,\phantom{a} t\ge t_0$.}

$\phantom{aaaaaaaaaaaaaaaaaaaaaaaaaaaaaaaaaaaaaaaaaaaaaaaaaaaaaaaaaaaaaaaaaaaa} \Box$

{\bf Theorem 3.3}. {\it Let $\phi_1(t)$ be  a solution of Eq.  (2.8) on the interval  $[t_0;+\infty)$, and  $\phi_0(t)$ be a solution of Eq. (1.1) such, that

\noindent
А$_1) \phantom{a} \phi_1(t_0) \ge \phi_0(t_0) > 0,\phantom{a} \phi'_1(t_0) > \phi_0(t_0) \ge 0$  or

$\phi_1(t_0) \le \phi_0(t_0) < 0, \phantom{a}\phi'_1(t_0) < \phi_0(t_0) \le 0$;

\noindent
B$_1) \phantom{a} p_0(t_0;\phi_0(t_0))\frac{\phi'_0(t_0)}{\phi_0(t_0)} < p_1(t_0;\phi_1(t_0))\frac{\phi'_1(t_0)}{\phi_1(t_0)}$.

\noindent
Let in addition the following conditions be satisfied:

\noindent
C$_1)  \phantom{a} p_0(t;w)\equiv p_1(t;w)$ is a  non increasing by  $w$  on the interval $(-\infty;0]$  and non decreasing by $w$  on the interval $[0;+\infty)$ function;

\noindent
D$_1)  \phantom{a} r_1(t;w_1)\le r_0(t;w)\le0$ for  $t\ge t_0, \phantom{a}|w|\le|w_1|,\phantom{a} w,\phantom{a}w_1 \in (-\infty;+\infty)$;

\noindent
E$_1)  \phantom{a} \frac{q_0(t;w)}{p_0(t;w)}\le \frac{q_1(t;w_1)}{p_1(t;w_1)}$ for  $t\ge t_0, \phantom{a}|w|\le|w_1|,\phantom{a} w,\phantom{a} w_1 \in (-\infty;+\infty)$.

\noindent
Then  $\phi_0(t)$  exists on the interval $[t_0;+\infty)$, and the function  $|\phi_0(t)|$ is positive and non decreasing.
}{

Proof. Let $[t_0;T)$ be the maximum existence interval  for  $\phi_0(t)$.  Show, that
$$
\phi_0(t)\ne 0,\phantom{aaa} t\in [t_0;T).  \eqno (3.9)
$$
Suppose  this relation is false. Then it follows from  А$_1)$  that for some  $T_1\in(t_0;T)$
$$
\phi_0(t)\ne 0,\phantom{aaa} t\in [t_0;T_1); \eqno (3.10)
$$
$$
\phi_0(T_1) = 0.    \eqno (3.11)
$$
It follows from  (3.10)  that $y_0(t)\equiv p_0(t;\phi_0(t))\frac{\phi'_0(t)}{\phi_0(t)}$
exists at least on the interval $[t_0;T_1)$. It follows from  A$_1)$, that
$y_0(t) \ge 0,\phantom{a} t\in [t_0;T_1)$. By  (2.2) from here and A$_1)$ it follows  $\phi_0(T_1)\ne0$,  which contradicts (3.11). The obtained contradiction proves (3.9).
It follows from (3.9) that $y_0(t)$  exists on the interval $[t_0;T)$.
Then since  $y_0(t_0)\ge 0$  by  (2.5) it follows from  D$_1)$ that
$$
y_0(t) \ge 0,\phantom{aaa} t\in [t_0;T).  \eqno (3.12)
$$
Since  by condition of the theorem $\phi_1(t)\ne 0, \phantom{a}t\ge t_0$,   $y_1(t)\equiv p_1(t;\phi_1(t))\frac{\phi'_1(t)}{\phi_1(t)}$ exists on the interval $[t_0;+\infty)$.
It follows from  B$_1)$, that $y_1(t_0) > y_0(t_0)$. On the strength of  Lemma 2.5 it follows from here, from A$_1)$  and C$_1)$  - E$_1)$, that
$$
y_0(t) < y_1(t),\phantom{aaa} t\in [t_0;T).     \eqno (3.13)
$$
Let us show, that $T=+\infty$. Suppose $T<+\infty$.Then from В$_1)$  and  (3.13) it follows:
$$
\int\limits_{t_0}^T\frac{y_0(\tau)}{p_0(\tau;\phi_0(\tau))} d\tau \le \int\limits_{t_0}^T\frac{y_1(\tau)}{p_0(\tau;0)} d\tau < +\infty.
$$
Using Lemma 2.2 from here we conclude, that $[t_0;T)$
is not the maximum existence interval for $y_0(t)$.
But in the other hand since  $[t_0;T)$  is the maximum  existence interval for  $\phi_0(t)$ the set   $[t_0;T)$  is the maximum  existence interval for  $y_0(t)$.
We came to the contradiction. The obtained contradiction shows  that  $T=+\infty$. Thus  $\phi_0(t)$  exists on the interval $[t_0;+\infty)$.
Due to  (2.2) it follows from  А$_1)$  and  (3.12)
that the function  $|\phi_0(t)|$ is positive and nondecreasing. The theorem is proved.

{\bf Definition 3.1}. {\it  A solution  $\phi_0(t)$  of Eq.  (1.1) is called
singular oscillatory of second kind, if the existence domain of the function $\phi_0(t)$ is a bounded set, and if $\phi_0(t)$ infinitely many times changes  sign}.

{\bf Theorem 3.4}. {\it Let the following conditions be satisfied:

\noindent
A$_2)\phantom{a} p_0(t;w)\ge P(t),\phantom{a} \frac{q_0(t;w)}{p_0(t;w)} \ge Q(t), \phantom{a}t\ge t_0, \phantom{a}w\in (-\infty;+\infty)$;

\noindent
B$_2)\phantom{a} r_0(t;w) \ge 0, \phantom{a}t\ge t_0,\phantom{a} w\in (-\infty;+\infty)$,

\noindent
Then for each  $\phi_{(0)}$  and  $\phi_{(1)}$  a non-extendable  on the interval $[t_0;+\infty)$
solution $\phi_0(t)$  of Eq. (1.1),
satisfying the initial value conditions $\phi_0(t_0)=\phi_{(0)},\phantom{a}  \phi'_0(t_0)=\phi_{(1)}$, is singular oscillatory of second kind.
}

Proof. Let $[t_0;T)\phantom{a} (T < +\infty)$ be  the maximum  existence interval for  $\phi_0(t)$. Then it is evident that
$$
\sup\limits_{\xi\in [T_1;T)}|\phi_0(\xi)| > 0,\phantom{a} T_1 \in [t_0;T). \eqno (3.14)
$$
(otherwise $\phi_0(t)$ will be extended by zero, i.e. $\phi_0(t)\equiv 0$)
Let us show, that for each $T_1 \in [t_0;T)$ the function $\phi_0(t)$ has a zero  on the  interval $[T_1;T)$.
Suppose, that for some  $T_0 \in [t_0;T)$ the function  $\phi_0(t)$  has no zero on the interval $[T_0;T)$.
 Suppose then  $\phi_0(t) >~ 0, \phantom{a} t \in  [T_0;T)$   (the proof in the case  $\phi_0(t) < 0,\phantom{a} t\in [T_0;T)$, by analogy). By (2.6) it follows  from here,  from A$_2)$ and  B$_2)$ that
$$
\phi_0(t) \le \frac{p_0(T_0;\phi(T_0))}{P(t)}|\phi'(T_0)|\exp\biggl\{-\int\limits_{T_0}^t Q(\tau) d\tau\biggr\},\phantom{aaa} t\in [T_0;T).
$$
Therefore, $\phi_0(t)$ is bounded. Then due to (2.6) and  (2.7)
there exists  finite limits  $\lim\limits_{t\to T - 0}\phi_0(t)$,  $\lim\limits_{t\to T - 0}\phi'_0(t)$.                 It follows from here  that  $[t_0;T)$ is not the maximum  existence interval for  $\phi_0(t)$. The obtained contradiction shows  that for every  $T_1 \in [t_0;T)$ the function  $\phi_0(t)$ has a zero on the interval $[T_1;T)$. On the basis of Lemma 2.7 we conclude that  from here, from B$_2)$ and (3.14) it follows that
 $\phi_0(t)$  is a singular oscillatory solution of second kind. The theorem is proved.

{\bf Definition 3.2}. {\it  A solution of Eq.  (1.1) is called oscillatory, if it exists on the interval $[t_0;+\infty)$ and in every  neighborhood  of  $+\infty$ changes  sign.}

{\bf Definition 3.3}. {\it  A solution $\phi(t)$ of Eq. (1.1) is called singular oscillatory of first kind, if it exists on the interval $[t_0;+\infty)$, $supp\hskip 3pt \phi(t)$ is bounded and $\phi(t)$  infinitely many times changes sign.}

Let for every $\varepsilon > 0$
the functions  $p_\varepsilon(t), \phantom{a} q_\varepsilon(t)$   and   $r_\varepsilon(t)$  be  real valued  and continuous on the interval $[t_0;+\infty)$, and let
$p_\varepsilon(t)> 0,\phantom{a} t\ge t_0, \phantom{a}\varepsilon > 0$. Consider the family of equations.
$$
(p_\varepsilon(t)\phi'(t))' +q_\varepsilon(t)\phi'(t) + r_\varepsilon(t)\phi(t) = 0,\phantom{aaa}
t\ge t_0,\phantom{aaa} \varepsilon > 0. \eqno (3.15_\varepsilon)
$$

{\bf Theorem 3.5}. {\it Let the following conditions hold:

\noindent
A$_3)\phantom{a} r(t;w) \ge 0,\phantom{a} t\ge t_0,\phantom{a} w\in (-\infty;+\infty)$;

\noindent
B$_3)$    there exists  $\varepsilon_0 > 0$ such,
that for every  $\varepsilon \in (0;\varepsilon_0]$
$$
p_0(t;w) \le  p_\varepsilon(t),\phantom{a}
 \frac{q_0(t;w)} {p_0(t;w)}\le \frac{ q_\varepsilon(t)}{p_\varepsilon(t)},\phantom{a}
 r_0(t;w) \ge r_\varepsilon(t), \phantom{a}\mbox{for}\phantom{a} |w| \ge \varepsilon ,\phantom{a} t\ge t_0,
$$
and  Eq.   $(3.15_\varepsilon)$  is oscillatory;

\noindent
C$_3)$  there exists $N > 0$   such, that

\noindent
C$^1_3)$ \phantom{a}$p_0(t;w) \le P(t), \phantom{a}\frac{q_0(t;w)}{p_0(t;w)}\le Q(t)$ for   $|w|\le N, \phantom{a}t\ge t_0$ and
$$
\int\limits_{t_0}^{+\infty}\exp\biggl\{- \int\limits_{t_0}^\tau Q(s)d s\biggr\}\frac{d\tau}{P(\tau)} = +\infty; \eqno (3.16)
$$

\noindent
C$^2_3)$   for every  $\varepsilon \ge N$  the following inequalities hold:   $p_0(t;w) \le p_\varepsilon(t),\phantom{a} \frac{q_0(t;w)}{p_0(t;w)}\le q_\varepsilon(t),\\ \phantom{a} r_0(t;w) \ge r_\varepsilon(t)$ for  $N \le |w| \le \varepsilon$ and
$$
\int\limits_{t_0}^{+\infty}\frac{d\tau}{P(\tau)}\int\limits_{t_0}^\tau\exp\biggl\{-\int\limits_s^\tau q_\varepsilon(\xi)d\xi\biggr\}r_\varepsilon(s)ds = +\infty. \eqno (3.17)
$$
Then each  existing on the interval $[t_0;+\infty)$ nontrivial solution of Eq.  (1.1)  either oscillatory, or singular oscillatory of first kind.
}

Proof. Let $\phi_0(t)$ be a solution of Eq. (1.1) such that $supp \hskip 3pt  \phi_0(t)$ is unbounded
on the interval $[t_0;+\infty)$. Let us show that $\phi_0(t)$  has arbitrarily large zeros.
Suppose, that it is not so, i. e.  there exists  $t_1\ge t_0$ such that  $\phi_0(t) \ne0,\phantom{a} t \ge t_1$. Suppose then  $\phi_0(t)> 0, \phantom{a} t\ge t_1$     (the proof in the case  $\phi_0(t)< 0, \phantom{a} t\ge t_1$   by analogy). Due to (2.6) it follows from  A$_3)$, that there can be one of the following three cases.

\noindent
$\alpha)\phantom{a} \phi'_0(t) \ge 0, \phantom{aaa}t\ge t_1$;

\noindent
$\beta)$ there exists $t_2\ge t_1$ such,
that  $\phi_0(t) \le N,\phantom{a} t\ge t_2,\phantom{a} \phi'_0(t_2) < 0$;

\noindent
$\gamma)$ there exists  $t_2\ge t_1$ such,
that  $\phi_0(t) \ge N, \phantom{a}t\ge t_2, \phantom{a}\phi'_0(t_2) < 0$.

\noindent
Let the case  $\alpha)$ be satisfied and let  $0 < \varepsilon < \min \{\phi_0(t_1); \varepsilon_0\}$.                     Then  $y_0(t)\equiv p_0(t;\phi_0(t))\frac{\phi'_0(t)}{\phi_0(t)}$  is a solution  of Eq.  (2.1) on the interval $[t_2;+\infty)$,    and  $y_0(t) \ge 0, \phantom{a}t\ge t_1$. By virtue of  Lemma~2.6 it follows from here and  from B$_3)$, that the Riccati's equation
$$
y'(t) + \frac{y^2(t)}{p_\varepsilon(t)} +\frac{q_\varepsilon(t)}{p_\varepsilon(t)}y(t) +r_\varepsilon(t)=0,\phantom{aaa} t\ge t_1,
$$
has a solution on the interval $[t_1;+\infty)$. Consequently, corresponding equation  $(3.15_\varepsilon)$ is not oscillatory, which contradicts  C$_3)$. The obtained contradiction shows, that  $\phi_0(t)$ has arbitrary large zeroes. Let the condition $\beta)$ holds. Then it follows from  C$^1_3)$, that
$$
\int\limits_{t_2}^t\exp\biggl\{ -\int\limits_{t_2}^\tau
\frac{q_0(s;\phi_0(s))}{p_0(s;\phi_0(s))}\biggr\}
\frac{d\tau}{p_0(\tau;\phi_0(\tau))} \ge
 \int\limits_{t_2}^t\exp\biggl\{ -\int\limits_{t_2}^\tau Q(s) d s\biggr\}\frac{d\tau}{P(\tau)},\phantom{aaa} t\ge t_2 \eqno (3.18)
$$
It is easy to show that from (3.16) it follows equality
$
\int\limits_{t_2}^{+\infty}\exp\biggl\{- \int\limits_{t_2}^\tau Q(s) d s\biggr\}\frac{d\tau}{P(\tau)} = +\infty.
$
By  (2.7) from here, from A$_3)$  and  (3.18) it follows: $\lim\limits_{t\to +\infty}\phi_0(t) = - \infty$,  which contradicts the supposition: $\phi_0(t) > 0, \phantom{a}t\ge t_1$.
The obtained contradiction shows  that $\phi_0(t)$  has arbitrary large zeroes.
 Let the case $\gamma)$ be satisfied.  Then from C$^2_3)$ it follows:
$$
\int\limits_{t_2}^t \frac{d\tau}{p_0(\tau;\phi_0(\tau))}\int\limits_{t_2}^\tau\exp\biggl\{- \int\limits_s^\tau\frac{q_0(\xi;\phi_0(\xi))}{p_0(\xi;\phi_0(\xi))}d\xi\biggr\}r_0(s;\phi_0(s)) \phi_0(s) ds \ge\phantom{aaaaaaaaaaaaaaaaaaaaaaaa}
$$
$$
\phantom{aaaaaaaaaaaaaaaaaaaaaaaaaaaaaaaaaaaaaaaaa} \ge\int\limits_{t_2}^t \frac{d\tau}{p_\varepsilon(\tau)}\int\limits_{t_2}^\tau\exp\biggl\{- \int\limits_s^\tau q_\varepsilon(\xi) \biggr\} r_\varepsilon(\tau) d \tau
$$
By  (2.7) from here  and from(3.17) it follows: $\lim\limits_{t\to +\infty} \phi_0(t) = -\infty$, which contradicts the supposition $\phi_0(t) > 0,\phantom{a} t\ge t_1$. The obtained contradiction shows, that $\phi_0(t)$  has arbitrary large zeroes. Thus we showed  that a solution  of Eq. (1.1),
existing on the interval $[t_0;+\infty)$ and  $supp$  of which is an unbounded set, has arbitrary large zeroes.
Due to  Lemma 2.7 it follows  from here  and from A$_3)$  that  $\phi_0(t)$ is oscillatory.
Let  $\phi_0(t)$ be a  solution of Eq. (1.1) with the bounded support on the interval $[t_0;+\infty)$, and let  $\phi_0(t) = 0,\phantom{a} t\ge T$;
$$
\max\limits_{\xi\in[t_1;T]}|\phi_0(\xi)| > 0, \phantom{a}t_1 \in [t_0; T).   \eqno (3.19)
$$
Show  that for each   $T_1\in [t_0;T)$ the function $\phi_0(t)$ has a zero  on the interval $[T_1;T)$. Suppose  it is not so, i .e. there exists  $T_0\in [t_0;T)$ such, that  $\phi_0(t)$  does not vanish on the interval $[T_0;T)$. Since $\phi_0(T)=0$   on the strength of Lagrange's mean value theorem  there exists  $\xi\in [T_0;T)$  such  that  $\phi'_0(\xi)\ne 0$, and  $sign \hskip 3pt \phi_0(T_0) = - sign \hskip 3pt \phi'_0(\xi)$.
By  (2.6) it follows from here  and from A$_3)$, that $\phi'_0(T)\ne 0$. But, in the other hand,
since   $\phi_0(t)=0$   for  $t\ge T$ we have  $\phi'_0(T) = 0$.
We came to the contradiction. Consequently for every  $T_1\in [t_0;T)$ the function  $\phi_0(t)$ has a zero on the interval $[T_1;T)$.
Due to Lemma 2.7  it follows from here, from  A$_3)$ and (3.19), that $\phi_0(t)$ is a singular oscillatory solution of first kind for  Eq. (1.1).  The  theorem is proved.

{\bf Remark 3.2}. {\it If the solution $\phi_0(t)$  of Eq.  (1.1),
satisfying the initial value conditions $\phi_0(t_0) = \phi'_0(t_0)=0 \phantom{a}$ is unique,  $(\phi_0(t)\equiv0)$,
in particular, if  $p_0(t;w),\phantom{a} q_0(t;w)$  and  $r_0(t;w)$  satisfy the conditions of the remark 1.1, then eq. (1.1) has no singular oscillatory solutions of first kind.}

{\bf Theorem 3.6}. {\it Let the following conditions  be satisfied:

\noindent
A$_4)\phantom{a} r_0(t;w) \ge0$  for  $t\ge t_0, \phantom{a}w\in (-\infty;+\infty)$;

\noindent
B$_4)\phantom{a} p_0(t;w), \frac{q_0(t;w)}{p_0(t;w)},\phantom{a} - r_0(t;w)$  are non increasing by  $w$  on the  interval $(-\infty;0]$ and nondecreasing by  $w$  on the interval $[0;+\infty)$ functions.

\noindent
Then for every $\phi_{(0)}$ and  $\phi_{(1)}$ the solution  $\phi_0(t)$ of Eq. (1.1),
satisfying the initial  conditions:
$$
\phi_0(t_0) = \phi_{(0)},\phantom{a}  \phi'_0(t_0) = \phi_{(1)}, \phantom{a}  \eqno (3.20)
$$
exists on the interval $[t_0;+\infty)$.}

Proof. Let  $\phi_0(t)$ be a  solution of Eq.  (1.1),
satisfying the initial  conditions (3.19),and let  $[t_0;T)$ be the maximum interval of existence for $\phi_0(t)$.  We should show, that  $T=+\infty$.
Suppose  $T<+\infty$. Two cases are possible:

\noindent
$\alpha)$ there exists  $t_1\in [t_0;T)$ such, that  $\phi_0(t) \ne 0,\phantom{a} t\in [t_1;T)$;

\noindent
$\beta)$  there exists infinite sequence  $t_0< t_0 < t_2 < ... < T$ such, that $\lim\limits_{k\to +\infty}t_k = T, \phantom{a}\phi_0(t_k) =~ 0,\linebreak k= 1,2, ...$.

\noindent
 Suppose the case  $\alpha)$ holds,  and   $\phi_0(t) > 0,\phantom{a} t\in [t_1; T)$  (the proof in the case  $\phi_0(t) < 0, \linebreak t\in~ [t_1; T)$,  by analogy). If  $\phi'_0(t)\le 0$, then by virtue of (2.7) it follows from  A$_4)$, that $\phi_0(t)$  has finite (nonnegative) limit when $t \to T -0$. Then by virtue of  (2.6) the function  $\phi'_0(t)$
also has finite limit when  $t \to T -0$. Therefore, $\phi_0(t)$
is continuable to the right at  $T,$  so  $[t_0; T)$
is not the maximum  existence  interval for $\phi_0(t)$.
The obtained   contradiction shows that  $T=+\infty$. Suppose  $\phi'_0(t_1) > 0$.
Two subcases are possible:

\noindent
$\alpha_1)\phantom{a}  \phi'_0(t) \ge 0,\phantom{a} t\in [t_1;T)$;

\noindent
$\beta_1)\phantom{a} \phi'_0(t) \ge 0$,\phantom{a} for  $t\in [t_1;T_1],\phantom{a} \phi'_0(t) < 0$, for  $\phantom{a} t\in (T_1;T)$, for some  $T_1 \in (t_1; T)$ (by virtue of (2,2) and (2.5) it follows from A$_4)$, that if $\phi_0'(\widetilde{t}_1) < 0$ for some $\widetilde{t}_1 \in (t_1;T]$, then $\phi_0'(t) < 0, \phantom{a} t\in [\widetilde{t}_1;T]$).

\noindent

In the case  $\alpha_1)$ the function  $\phi_0(t)$ is nondecreasing  on the interval  $(t_1; T)$.
Then  (2.6) it follows from A$_4)$  and B$_4)$, that
$$
\phi'_0(t) \le \phi'_0(t_1) \exp\biggl\{ - \int\limits_{t_1}^t\frac{q_0(\tau;\phi_0(t_1))}{p_0(\tau;\phi_0(t_1))}d\tau\biggr\}, \phantom{aaa}t\in [t_1; T).
$$
Therefore,  $\phi_0(t)$ has a finite limit when  $t\to T-0$.
It is evident that the same we have  in the case  $\beta_1)$. Then by  (2.6) $\phi'_0(t)$  has a finite limit  when  $t\to T-0$.  So, $[t_0; T)$ is not the maximum  existence interval for $\phi_0(t)$.
The obtained contradiction shows  that $T=+\infty$.
 Assume the case $\beta)$ takes place. If  $\phi_0(t)$ is bounded, then, it is evident, that the functions $\frac{1}{p_0(t,\phi_0(t))}, \phantom{a} \frac{q_0(t,\phi_0(t))}{p_0(t,\phi_0(t))}, \phantom{a} r_0(t,\phi_0(t))$
are bounded in the interval $[t_0;T)$.
By  (2.7) it follows from here, that $\phi_0(t)$ has finite limit when $t\to T-0$.
 Then arguing similarly to the above (when we are dealing with the  cases $\alpha_1)$ and $\beta_1)$) we conclude, that $T=+\infty$.  Assume  $\phi_0(t)$ is not bounded. Two subcases are possible:

\noindent
$\alpha_2)\phantom{a} \stackrel{ \overline{\lim}}{_{_{t\to T - 0}}} \phi_0(t) = +\infty$;

\noindent
$\beta_2)\phantom{a} \stackrel{ \underline{\lim}} {_{_{t\to T - 0}}} \phi_0(t) = -\infty$.

\noindent
Let the case $\alpha_2)$ be satisfied   (the proof in the case  $\beta_2)$ by analogy).
Let $\phi_1(t)$  be  the solution of Eq. (1.1) with  $\phi_1(T) =1, \phantom{a} \phi'_1(T) =-1$.
Then there exists $\zeta \in [t_1;T)$ such  that $\phi_1(t)$ exists on the interval $[\zeta;T]$, and
$$
\phi_1(t) > 0,\phantom{a} \phi'_1(t) < 0,\phantom{aaa} t\in [\zeta;T].   \eqno (3.21)
$$
It follows from  $\beta)$  and $\alpha_2)$ that there exists a point  $\zeta_1 (\in [\zeta;T))$
of local maximum for  $\phi_0(t)$  such  that  $\phi_0(\zeta_1) > \max\limits_{\xi \in [\zeta;T]}\phi_1(\xi)$    and  $\zeta < t_m < \zeta_1$ for some $m$.
Since $t_m < \zeta_1 < \\ < t_n$ for some $n > m$  and  $\phi_0(t_m) = \phi_0(t_n)=0< \min\{\phi_1(t_m), \phantom{a}\phi_1(t_n)\}$, but  $\phi_0(\zeta_1) >\linebreak > \max\{\phi_1(t_m), \phantom{a}\phi_1(t_n)\}$,
there exist $\xi_1$ and $\xi_2$ such that   $\phi_0(\xi_k) =  \phi_1(\xi_k),\phantom{a} k=1,2$
and
$$
\phi_0(t) > \phi_1(t), \phantom{aaa}t\in (\xi_1;\xi_2). \eqno (3.22)
$$
(see pict. 1).
Let  $y_j(t)\equiv  p_0(t;\phi_j(t))\frac{\phi'_j(t)}{\phi_j(t)},\phantom{a} j=0,1$.
Then by virtue of $(2.12_1)$ the following equality holds:
$$
y_1(t) - y_0(t) = [y_1(\zeta_1) - y_0(\zeta_1)]\exp\biggl\{- \int\limits_{\zeta_1}^t\frac{y_0(\tau) + y_1(\tau) + q_0(\tau;\phi_0(\tau))}{p_0(\tau;\phi_0(\tau))}d\tau\biggr\} -\phantom{aaaaaaaaaaaaaa}
$$
$$
- \int\limits_{\zeta_1}^t\exp\biggl\{-\int\limits_{\tau}^t\frac{y_0(s) + y_1(s) + q_0(s;\phi_0(s))}{p_0(s;\phi_0(s))}ds\biggr\}\biggl[\biggl(\frac{1}{p_1(\tau;\phi_1(\tau))} - \frac{1}{p_0(\tau;\phi_0(\tau))}\biggr) y_1^2(\tau) + \biggr.
$$
$$
\phantom{aaaaa}\biggl. +\biggl(\frac{q_1(\tau;\phi_1(\tau))}{p_1(\tau;\phi_1(\tau))} - \frac{q_0(\tau;\phi_0(\tau))}{p_0(\tau;\phi_0(\tau))}\biggr) y_1(\tau) + r_1(\tau;\phi_1(\tau)) - r_0(\tau;\phi_0(\tau))\biggr]d\tau = 0. \eqno (3.23)
$$

\vskip 40pt

\begin{picture}(120,100)
%\put(145,0){\vector(0,1){70}}
\put(40,30){\vector(1,0){330}}
%\put(174,15){\bf \line(1,0){180}}
\put(58,30){\circle*{3}}
%\put(145,5){\circle*{3}}
\put(290,30){\circle{3}}
\put(288,17){$T$}
\put(58,17){$t_0$}
\put(360,17){$t$}
\put(165,0){$\small \mbox{pict. 1}$}
%\put(115,61){\bf \vector(-1,0){8}}
\multiput(180,85)(1,-.25){20}{.}
\multiput(200,80)(1,-.2){20}{.}
\multiput(218,76)(1,-.3){20}{.}
\multiput(238,70)(1,-.4){20}{.}
\multiput(258,62)(1,-.6){20}{.}
\multiput(278,50)(1,-1){15}{.}
\multiput(292,35)(1,-.6){20}{.}
\multiput(288,30)(0,4){20}{.}
\multiput(200,30)(1,.7){5}{.}
\multiput(204,33)(.5,.7){10}{.}
\multiput(209,40)(.5,1){10}{.}
\multiput(214,50)(.5,1.4){10}{.}
\multiput(219,64)(.25,1.1){20}{.}
\multiput(224,85)(.25,1){15}{.}
\multiput(228,100)(.25,.5){10}{.}
\multiput(231,105)(.25,.3){10}{.}
\multiput(233,108)(.25,.1){10}{.}
\multiput(235,109)(.25,.0){10}{.}
\multiput(238,108)(.25,-.1){10}{.}
\multiput(240,106)(.25,-.3){16}{.}
\multiput(242,104)(.25,-.5){20}{.}
\multiput(247,93.5)(.25,-.7){26}{.}
\multiput(254,75)(.25,-.9){26}{.}
\multiput(259,57)(.25,-.8){26}{.}
\multiput(265.5,36)(.25,-.9){10}{.}
\multiput(222,30)(0,4){14}{.}
\multiput(236,30)(0,4){20}{.}
\multiput(257.5,30)(0,4){10}{.}
\put(224,30){\circle*{3}}
\put(238,30){\circle*{3}}
\put(259,30){\circle*{3}}
\put(220,17){$\xi_1$}
\put(234,17){$\zeta_1$}
\put(255,17){$\xi_2$}
\put(188,47){$\phi_0(t)$}
\put(265,63){$\phi_1(t)$}
\put(6,-20){$\small{\small \mbox{An illustration to the part of proof of Theorem 3.6, related to the case}\hskip 4 pt \alpha_2). \hskip 4pt \mbox{Here}}  \hskip 4pt  [t_0;T) $}
\put(6,-32){$\small  \mbox{is the maximum existence interval for} \hskip 4pt \phi_0(t).\hskip 4pt \mbox{The "cap" \hskip 4pt in this picture is a part of the} $}
\put(6,-44){$\small  \mbox{graph of}\hskip 4pt \phi_0(t). \mbox{The decreasing curve is a part of the graph of}\hskip 4pt \phi_1(t);\hskip 4pt  \xi_1 \hskip 4pt  \mbox{and} \hskip 4pt \xi_2 \hskip 4pt  \mbox{are} $}
\put(6,-56){$\small  \mbox{points of intersection of graph of}\hskip 4pt \hskip 4pt \phi_0(t)  \hskip 4pt  \mbox{and} \hskip 4pt \phi_1;   \hskip 4pt   \zeta_1  \hskip 4pt  \mbox{is a local maximum point of} \hskip 4pt \phi_0(t).$}
\end{picture}

\vskip 80pt
\noindent
It is evident that
$$
y_1(\zeta_1) < y_0(\zeta_1)\phantom{a} (=0). \eqno (3.24)
$$
Since $\phi_0(t) \ge \phi_1(t) > 0$   on the interval $[\zeta_1;\xi_2]$    it follows from the conditions of the theorem that
$$
p_0(t;\phi_0(t)) \ge p_0(t;\phi_1(t)), \phantom{a} r_0(t;\phi_0(t)) \le r_0(t;\phi_1(t)),\phantom{a} \frac{q_0(t;\phi_0(t))}{p_0(t;\phi_0(t))}\ge \frac{q_0(t;\phi_1(t))}{p_0(t;\phi_1(t))},\phantom{a} t\in [\zeta_1;\xi_2].
$$
Then since  $y_1(t) < 0$ on the interval $[\zeta_1;\xi_2]$  it follows  from (3.23) and (3.24),
that $y_1(\xi_2) < y_0(\xi_2)$.  Then  $\phi'_1(\xi_2) < \phi'_0(\xi_2)$.
It follows from here that  $\phi'_1(t) < \phi'_0(t), \phantom{a}t\in[\xi_3;\xi_2]$, for some  $\xi_3\in (\zeta_1; \xi_2)$. Therefore  $\int\limits_{\xi_3}^{\xi_2}(\phi_1(\tau) - \phi_0(\tau))'d\tau < 0$
or, which is the same, $\phi_1(\xi_3) >\\> \phi_0(\xi_3)$ (since $\phi_1(\xi_2) = \phi_0(\xi_2),\phantom{a} see \hskip 4pt pict. 1$)  which contradicts  (3.22).
The obtained contradiction shows, that $T=+\infty$.  The theorem is proved.

\vskip 12pt
{\bf \centerline { \S 4. Some applications}}
\vskip 12pt

Consider the Emden - Fowler equation  (see  [2], p. 171):
$$
(t^\rho\phi'(t))' - t^\sigma \phi^n(t) = 0, \phantom{a}t\ge t_0 > 0,\phantom{a} n > 1, \phantom{a}\rho,\phantom{a} \sigma \in (-\infty;+\infty). \eqno (4.1)
$$
Along with this equation consider the equation
$$
(t^\rho\phi'(t))' - t^\sigma |\phi(t)|^{n-1}\phi(t) = 0,\phantom{aaa} t\ge t_0. \eqno (4.2)
$$
Here  $p_0(t;w)\equiv t^\rho ,\phantom{a} q_0(t;w) \equiv 0,\phantom{a} r_0(t;w) \equiv - t^\sigma |w|^{n-1}$. Set: $P(t) \equiv t^\rho , \phantom{a}Q(t) \equiv \\ \equiv 0, \phantom{a}R(t) \equiv - t^\sigma$.                                 Assume $\rho >1$.   Then
$$
F(t_0; t; c_1; c_2) \le |c_1| \exp\biggl\{\frac{c_2}{\rho - 1}(t_0^{1-\rho} - t^{1-\rho}) +
 \frac{t^{\sigma + 2 -\rho} - t_0^{\sigma + 2 -\rho}}{(\sigma+1)(\sigma + 2 - \rho)} - \frac{t_0^{\sigma + 2 - \rho}}{(\sigma+1)(\rho - 1)}    \biggr\}.
$$
It follows from here that if $-1 < \sigma < \rho - 1$, then
$$
F(t_0; t; c_1; c_2) \le |c_1| \exp\biggl\{\frac{c_2}{\rho - 1}t_0^{1-\rho} - \frac{t_0^{\sigma + 2 -\rho}}{(\sigma+1)(\sigma + 2 - \rho)}\biggr\} \stackrel{def}=A(t_0;c_1;c_2),
$$
and if  $\sigma < -1$, then
$$
F(t_0; t; c_1; c_2) \le |c_1| \exp\biggl\{\frac{c_2}{\rho - 1}t_0^{1-\rho} - \frac{t_0^{\sigma + 2 -\rho}}{(\sigma+1)(\rho -1)}\biggr\} \stackrel{def}=B(t_0;c_1;c_2).
$$
It is not difficult to see, that if
$$
-1 < \sigma < \rho - 2, \phantom{aaa}A\biggl(t_0;\phi_{(0)};t_0^\rho\frac{\phi_{(1)}}{\phi_{(0)}}\biggr) < 1, \eqno (4.3)
$$
 or
$$
\sigma < -1,  \phantom{aaa} B\biggl(t_0;\phi_{(0)};t_0^\rho\frac{\phi_{(1)}}{\phi_{(0)}}\biggr) < 1, \eqno (4.4)
$$
then for Eq. (4.2) the conditions of Theorem 3.1 are satisfied.
Therefore, for every  $\phi_{(0)},\phantom{a} \phi_{(1)}$,
satisfying  the  conditions  $\phi_{(0)}\ne0,\phantom{a} \frac{ \phi_{(1)}}{\phi_{(0)}} \ge 0$
and one of the conditions (4.3) or (4.4), the solution  $\phi_0(t)$
of Eq. (4.2) with  $\phi_0(t_0)= \phi_{(0)}, \phantom{a} \phi'_0(t_0)= \phi_{(1)}$,
exists on the interval $[t_0;+\infty)$, and if takes place  (4.3), then
$$
|\phi_0(t)| \le A\biggl(t_0;\phi_{(0)};t_0^\rho\frac{\phi_{(1)}}{\phi_{(0)}}\biggr),\phantom{aaa} t\ge t_0, \eqno (4.5)
$$
and if takes place (4.4), then
$$
|\phi_0(t)| \le B\biggl(t_0;\phi_{(0)};t_0^\rho\frac{\phi_{(1)}}{\phi_{(0)}}\biggr),\phantom{aaa} t\ge t_0; \eqno (4.6)
$$
the function $|\phi_0(t)|$ is positive and nondecreasing.
It is evident, that if $\phi_{(0)}> 0$, then $\phi_0(t)$  is a solution of Eq. (4.1) too.

{\bf Remark 4.1}. {\it If $\phi(t)$ is a positive (negative) solution of Eq. (4.2) (and $n = \frac{n_1}{n_2}$, where $n_1$ and $n_2$ are odd), then $\phi(t)$  is a solution of Eq. (4.1) too. Note that equations (4.1) and (4.2) are not equivalent for all $n$, e. g. the equation $\phi''(t) - \frac{1}{t}\phi^2(t) = 0$ is not equivalent to the equation  $\phi''(t) - \frac{1}{t}|\phi(t)|\phi(t) = 0$,
since the function $\phi_0(t)\equiv - \frac{2}{t}$ is a solution of the last equation on the interval $[1;+\infty)$, but not of the first one.}

Assume  $\rho = 0,\phantom{aaa} \sigma + n + 1 < 0$. Then (see  [2], p. 173)
the function
$
\phi_B(t)\equiv
$
\linebreak
$
\equiv \left[\frac{(\sigma+2)(\sigma+n + 1)}{(n - 1)^2}\right]^{\frac{1}{n - 1}} t^{-\frac{\sigma +2}{n-1}}
$
is a solution of Eq. (4.2). It is not difficult to see that for  \\$p(t;w)=p_1(t;w)\equiv 1,\phantom{a}  q(t;w)=q_1(t;w)\equiv 0,\phantom{a} r(t;w)=r_1(t;w)= - t^\sigma |w|^{n-1},\phantom{a} t\ge t_0,\\  w\in (-\infty; +\infty),\phantom{a} \rho= 0,\phantom{a} \sigma + n +1 <0$
for equations (2.8)and (4.2) the conditions  C$_1)$  - \\ - E$_1)$ of  Theorem 3.3 are satisfied.
Therefore, the solution $\phi_0(t)$ of Eq. (4.2) with  $\phi_0(t_0)\ne\\ \ne 0,\phantom{a} \frac{\phi'_0(t_0)}{\phi_0(t_0)} \le \frac{\phi'_B(t_0)}{\phi_B(t_0)}$,
exists on the interval $[t_0;+\infty)$, and $|\phi_0(t)|$  is positive and nondecrea-\\sing.
It is evident that if  $\phi_0(t_0)> 0$  then  $\phi_0(t)$ is a solution of Eq. (4.1). For $\phi_0(t_0)< 0$ the function  $\phi_0(t)$  will be a solution of Eq. (4.1) if  $n=\frac{n_1}{n_2}$,    where  $n_1$  and  $n_2$  are odd.
If $\rho \ne1$ then by invertible transformation
$$
\left\{
\begin{array}{l}
s=(\rho - 1)^{-1} t ^{\rho - 1}, \phantom{a} \phi = (\rho - 1)^{(\rho - \sigma - 2)/[(\rho -1)(n - 1)]}\frac{\psi}{s}\phantom{a} for \phantom{a} \rho > 1;\\
\phantom{a}\\
s = (1 - \rho)^{-1} t^{1 - \rho}, \phantom{a} \phi = (1 - \rho)^{-(\sigma + \rho)/[(n - 1)(1 - \rho)]} \psi, \phantom{a} for \phantom{a} \rho < 1.
\end{array}
\right.
$$
Eq. (4.2)  converts in (see [2], pp. 171, 172)
$$
\psi''(s) - s^{\sigma_1}|\psi (s)|^{n-1}\psi(s) = 0,\phantom{a} s\ge s_0 >0, \eqno (4.7)
$$
where
$$
\sigma_1=\left\{
\begin{array}{l}
\frac{\sigma + \rho}{\rho - 1} - (n + 3),\phantom{a} \rho > 1;\\
\phantom{a}\\
\frac{\sigma + \rho}{1 -\rho}, \phantom{a}\rho < 1.
\end{array}
\right.  \eqno (4.8)
$$
This transformation establishes a one-to-one correspondence  between solutions of equations (4.2) and  (4.7), and existing on the interval $[t_0;+\infty)$ solutions of Eq.  (4.2) transform in solutions of Eq. (4.7), existing on the interval $[s_0;+\infty)$. By  already  proved above for $\sigma_1 + n +1 < 0$  Eq.  (4.7) has two - parameter family of solutions on the interval $[s_0;+\infty)$. Therefore from (4.8) it  follows, that in the cases  $\rho > \max \{1, \sigma + 2\}$  and   $\frac{\sigma - 1}{n} + 1 < \rho  < 1$  Eq.  (4.2) has two - parameter family of solutions on the interval $[t_0;+\infty)$.

{\bf Definition 4.1}. {\it A solution  $(\phi_0(t), \psi_0(t))$  of Eq.  (1.2) is called
conditionally stable for  $t \to + \infty$,
if there exists one dimensional manifold  $S\ni (\phi_0(t_0), \psi_0(t_0))$ such,
that for every  $\varepsilon > 0$   and for every solution  $(\phi(t), \psi(t))$  of the system (1.2) there exists  $\delta > 0$  such,
that  $|\phi(t) - \phi_0(t)| + |\psi(t)  - \psi_0(t)| <\varepsilon$  for  $t\ge t_0$,
as soon as  $(\phi(t_0), \psi(t_0)) \in S$   and   $|\phi(t_0) - \phi_0(t_0)| + |\psi(t_0)  - \psi_0(t_0)| < \delta$  (see [10], p. 314).}

{\bf Definition 4.2}. {\it A solution $\phi_0(t)$  of Eq. (1.1) is called conditionally stable for $t \to  + \infty$,   if the corresponding solution  $(\phi_0(t), p_0(t;\phi_0(t))\phi'_0(t))$ of the system  (1.2)
is conditionally stable for $t \to + \infty$.}

Show that if $\rho > 1, \phantom{a}\sigma < -1$, then the solution  $\phi_0(t)\equiv 0$
of Eq. (4.1) is conditionally stable for  $t \to + \infty$.
Let  $S=\bigl\{(\phi_{(0)},  \phi_{(1)}) : 0 \le \phi_{(0)} \le  \exp\bigl\{\frac{t_0^{\sigma + 2 - \rho}}{(\sigma + 1)(\rho - 1)}\bigr\}, \phi_{(1)} =\\ = 0\bigr\}$, and let  $\phi(t)$ be a solution  of Eq.
(4.1) with  $(\phi(t_0), t_0^\rho \phi'_0(t_0)) \in S,\phantom{a} \phi(t_0) \ne 0$. Then  $B(t_0; \phi(t_0); 0) <1$.
By virtue of  Theorem 3.1 it follows from here  and from (4.4) that  $\phi(t)$  exists on the interval
$[t_0;+\infty)$ and satisfies the inequality
$$
|\phi(t)| \le \phi(t_0)\exp\biggl\{ - \frac{t_0^{\sigma + 2 - \rho}}{(\sigma + 1)(\rho - 1)}\biggr\},\phantom{aaa} t \ge t_0 > 0.   \eqno (4.9)
$$
By  (4.1) it follows  from her  that
$$
|t^\rho\phi'(t)| \le \phi(t_0)\exp\biggl\{ - \frac{t_0^{\sigma +
2 - \rho}}{(\sigma + 1)(\rho - 1)}\biggr\}\biggl(\frac{-t_0^{\sigma + 1}}{\sigma +1}\biggr),\phantom{aaa} t \ge t_0 > 0.   \eqno (4.10)
$$
Let $\varepsilon > 0$  be  fixed. Set:
$$
\delta = \delta(\varepsilon) \equiv\frac{\varepsilon}{2}\biggl(1 - \frac{t_0^{\sigma +1}}{\sigma +1}\biggr)^{-1}\exp\biggl\{  \frac{t_0^{\sigma +
2 - \rho}}{(\sigma + 1)(\rho - 1)}\biggr\}
$$
Let $|\phi(t_0)| + |t_0^\rho\phi'(t_0)| < \delta$.  Then it follows from (4.9) and (4.10) that
$$
|\phi(t)| + |t_0^\rho\phi'(t)| < \varepsilon \phantom{a}\mbox{for} \phantom{a}t\ge t_0.    \eqno (4.11)
$$
If  $\phi(t_0)=0$, then by virtue of  Remark 1.1 $\phi(t)\equiv 0$  and, therefore,
in this case the relation (4.11) also takes place.
Consequently, the solution  $(\phi_0(t), \psi_0(t))\equiv (0, 0)$ of the system
$$
\left\{
\begin{array}{l}
\phi'(t) =  \frac{\psi(t)}{t^\rho};\\
\phantom{a}\\
\psi'(t) = t^\sigma \phi^n(t)
\end{array}
\right.
$$
 is conditionally stable for $t\to \infty$. Then the solution  $\phi_0(t)\equiv 0$
of Eq. (4.1) is conditionally stable for  $t\to \infty$.
Taking into account  Remark 4.1 we summarize the obtained  result in the following form.

{\bf Theorem 4.1}. {\it The following assertions are valid.

\noindent
I).  Let $\rho > 1$, and let  $\phi_{(0)}$ and  $\phi_{(1)}$ satisfy the conditions: $\phi_{(0)}\ne0, \phantom{a} \frac{\phi_{(1)}}{\phi_{(0)}} \ge 0$
and one of the conditions (4.3), (4.4).
Then the solution $\phi_0(t)$  of Eq. (4.2),
satisfying the initial value conditions: $\phi_0(t_0) = \phi_{(0)},\phantom{a} \phi'_0(t_0) = \phi_{(1)}$,
exists on the interval $[t_0;+\infty)$.
The function $|\phi_0(t)|$ is positive and nondecreasing,
and if (4.3) holds, then the estimate (4.5) is valid,
and if (4.4) holds, then the estimate  (4.6) is valid.
If  $\phi_{(0)} > 0$  or if   $\phi_{(0)} < 0,\phantom{a} n = \frac{n_1}{n_2}$,
where  $n_1$  and $n_2$  are odd, then  $\phi_0(t)$ is a solution of  (4.1).

\noindent
II). Let  $\rho =0, \phantom{a}\sigma + n +1 < 0$. Then the solution  $\phi_0(t)$
of Eq.  (4.2) with   $\phi_0(t_0) \ne 0,\phantom{a} 0\le \\ \le \frac{\phi'_0(t_0)}{\phi_0(t_0)} < \frac{\phi'_B(t_0)}{\phi_B(t_0)}$, exists on the interval $[t_0;+\infty)$,
and  $|\phi_0(t)|$  is positive and nondecreasing.
If  $\phi_0(t_0) > 0$  or if  $\phi_0(t_0) <0,\phantom{a} n= \frac{n_1}{n_2}$,
where $n_1$  and  $n_2$  are odd, then  $\phi_0(t)$  is a solution of Eq.  (4.1).

\noindent
III). In the cases  $\rho > \max\{1, \sigma + 2\}$  and $\frac{\sigma - 1}{n} + 1 < \rho < 1$   Eq.  (4.1) has two - parameter family of solutions on the interval $[t_0; +\infty)$.

\noindent
IV).  If  $\rho >1$  and  $\sigma < -1$, then  $\phi_0(t)\equiv 0$
of eq. (4.1) is conditionally stable for  $t\to +\infty$.}

{\bf Remark 4.2}. In the case  $\rho = 0$ the existence of global solutions of Eq. (4.1), which are different by  their properties from described in assertion II of  theorem 4.1  (the Kneser's solutions) follows from  theorem 16.1 of book [1] (see [1], p. 371).

Let us compare  theorem 4.1 with the following result (see [11], p. 8).

{\bf Theorem*}. The following assertions hold:

\noindent
i). There exists $\varepsilon > 0$ such, that every solution of Eq. (4.7) with Cauchy initial conditions
$|\psi(s)| \le \varepsilon,\phantom{a} |\psi'(s| \le \varepsilon$ exists on the interval $[s_0;+\infty)$ if and only if $\sigma_1 < -n - 1$.

\noindent
ii). If $\sigma_1 \ge - n -1$, then every solution $\psi(s)$ of Eq. (4.7), satisfying $\psi(\tau)\psi'(\tau) > 0$ at some $\tau \ge s_0$ is non continuable on the interval $[s_0;+\infty)$.

In the assertions $I)$ and $II)$ of Theorem 4.1 the region of the initial values $\phi(t_0),\hskip 2pt \phi'(t_0)$ ($\psi(s_0), \phantom{a} \psi'(s_0)$) for which the solution $\phi(t)$ ($\psi(s)$) of Eq. (4.2) (of Eq. (4.7)) exists on the interval $[t_0;+\infty)$ ($[s_0;+\infty)$) is  describes by well - defined relationships, whereas from the assertion i) of  Theorem* we can not see exactly for which initial conditions $\psi(s_0), \phantom{a} \psi'(s_0)$ (except the trivial case $\psi(s_0) = \psi'(s_0) = 0$)
the solution $\psi(s)$  of Eq. (4.7) exists on the interval $[s_0;+\infty)$. From the assertion ii) of the theorem*  it follows, that in the assertion II) of  Theorem 4.1  the condition $\sigma + n +1 < 0$ can not be replaced by weaker condition $\sigma + n +1 \le 0$. In this sense  theorem 3.3 (which implies II)) is sharp. From the assertion ii) of  Theorem* and from (4.8) it follows, that in the assertion III) of  Theorem 4.1 the condition $\rho > \max\{1, \sigma + 2\}$  ($\frac{\sigma - 1}{n} +1 < 1$) can not be replaced by weaker condition $\rho \ge \max\{1, \sigma + 2\}$  ($\frac{\sigma - 1}{n} +1 \le 1$). In this sense  Theorem 3.1 (which is used in proof of III)) is sharp.

Let  $\lambda(t),\phantom{a} \mu(t)$  and $\nu(t)$  be continuous functions  on the interval $[t_0;+\infty)$
 and let   $\lambda(t) >\\> 0,\phantom{a} \mu(t) \ge 0,\phantom{a} \nu(t)\ge 0,\phantom{a} t\ge t_0$.
Consider the following Van der  Pol's  type equation (see [12]).
$$
(\lambda(t)\phi'(t))' + \mu(t)(\phi^2(t) - 1)\phi'(t) + \nu(t)\phi(t) = 0,\phantom{aaa} t\ge t_0. \eqno (4.12)
$$
Here  $p_0(t;w) \equiv \lambda(t),\phantom{a} q_0(t;w)\equiv \mu(t)(w^2 -1),\phantom{a} r_0(t;w)\equiv \nu(t)$
 satisfy all of the conditions of  Theorem 3.6.
Therefore, for each   $\phi_{(0)}$  and  $\phi_{(1)}$ the solution  $\phi_0(t)$ of Eq. (4.12), satisfying the initial value conditions: $\phi_0(t_0) = \phi_{(0)},\phantom{a} \phi'_0(t_0) = \phi_{(1)}$,
 exists on the interval $[t_0;+\infty)$. We put:
$$
P(t)=p_\varepsilon(t) \equiv \lambda(t),\phantom{a} Q(t)\equiv 0,\phantom{a} q_\varepsilon(t)=\mu(t)(\varepsilon^2 - 1),\phantom{a} r_\varepsilon(t)= \nu(t), \eqno (4.13)
$$
$t\ge t_0, \phantom{a} \varepsilon > 0, \phantom{a}N =1. $
It is not difficult to check, that if the following  conditions hold:
$$
 \mbox{a}^\circ). \int\limits_{t_0}^{+\infty}\frac{d\tau}{\lambda(\tau)} = \int\limits_{t_0}^{+\infty}\frac{d\tau}{\lambda(\tau)}\int\limits_{t_0}^\tau \exp\biggl\{- \int\limits_s^\tau\frac{\mu(\xi)(\varepsilon^2 - 1)}{\lambda(\xi)}\biggr\}\nu(s) d s = +\infty \phantom{a}\mbox{при}\phantom{a} \varepsilon \ge N;\phantom{aaaaaaaaaaaaaaaaaaaaaaaaaaaa}
$$

\noindent
$ \mbox{b}^\circ)$.  \phantom{a} for \phantom{a} $0< \varepsilon \le \varepsilon_0$
the equations
$$
(\lambda(t)\phi'(t))' + \mu(t)(\varepsilon^2 - 1)\phi'(t) + \nu(t)\phi(t) = 0,\phantom{aaa} t\ge t_0
$$
are oscillatory,
then for Eq. (4.12) with (4.13) the   conditions  A$_3)$   - B$_3)$ of  Theorem 3.5 hold.
Then  due to   Theorem  3.5 the solution  $\phi_0(t)$ either is oscillatory or is singular oscillatory of first kind. Since, it is evident, the functions
 $p_0(t;w) \equiv \lambda(t), \phantom{a}q_0(t;w)\equiv\\ \equiv \mu(t)(w^2 -1),\phantom{a} r_0(t;w)\equiv \nu(t)$  satisfy   the  conditions  of  Remark 1.1, the solution $\phi_0(t)$, satisfying the initial value conditions $\phi_0(t_0) = \phi_{(0)},\phantom{a} \phi'_0(t_0) = \phi_{(1)}$, is exactly one, and consequently, due to Remark 3.2 cannot be singular oscillatory of first kind. The obtained result we summarize in the following form.

{\bf Theorem 4.2}. {\it Let   $\lambda(t) > 0,\phantom{a} \mu(t)\ge 0,\phantom{a} \nu(t)\ge 0$ for  $t\ge t_0$.
Then for each   $\phi_{(0)}$   and  $\phi_{(1)}$  the solution  $\phi_0(t)$ of Eq.  (4.12), satisfying the initial value conditions: \\$\phi_0(t_0) = \phi_{(0)},\phantom{a} \phi'_0(t_0) = \phi_{(1)}$,
exists on the interval  $[t_0;+\infty)$. Moreover if in addition the conditions  $\mbox{a}^\circ)$ and  $ \mbox{b}^\circ)$ hold, then $\phi_0(t)$  is oscillatory.}

\vskip 22pt
{\bf \centerline {Conclusion}}
\vskip 12pt

We have used the Riccati equation method to investigate some classes of second order nonlinear ordinary differential equations. This method have made possible us to establish  four new global existence criteria for the mentioned classes of equations. We have proved two new oscillatory criteria for them as well. These criteria were used to the Emden - Fowler equation, having applications  in the astrophysics, and to the Wan der Pole type equation, which is applicable for studying the dynamics of dusty grain charge in dusty plasmas.

\vskip 22pt
{\bf \centerline {References}}
\vskip 12pt

1. Kiguradze I. T., Chanturia T. A. The asymptotic behavior of solutions of nonlinear\\ \phantom{aaaaa} ordinary differential equations, Moscow, ''Nauka'', 1990.

2. Bellman R. Stability theory of differential equations, Moscow, Izdatelstvo inostrannoj \phantom{aaaaa} literatury, 1954 (New York, Toronto, London,  McGRAW-HILL BOOK COMPANY,\linebreak \phantom{aaaaa} INC, 1053).

3. Hartman Ph. Ordinary differential equations. Moscow, ''Mir'', 1970  (JOHN WILLEY\linebreak \phantom{aaaaa}  and SONS, New York - London - Sydney 1064).

4. Reising R., Sansone G. Conty R. Qualitative theory of nonlinear differential equations, \phantom{aaaaa} Moscow, ''Nauka'', 1974 (EDIZIONI CREMONESE ROMA 1065).

5. Jordan D. W., Smith P. Nonlinear differential equations. Oxford University press,\\\phantom{aaaaa}  2007.

6. Cheng Y., Peng M., Zhang W. On the asymptotic behavior of the solutions of a \\ \phantom{aaaaa} class of second order nonlinear differential equations. Journal of Computational \\ \phantom{aaaaaa}and Applied
Mathematics 98 (1998) pp. 63 - 79.

7. Motohiko K., Kusano T. On a class of second order quasilinear  ordinary differen-\\ \phantom{aaaaa} tial equations. Hiroshima Math. J., v. 25, 1995, pp. 321 - 355.

8. Kusano T.,  Akio O. Existence and asymptotic behavior of positive solutions of \\ \phantom{aaaaa} second order quasilinear differential equations. Funkcialax Ekvacioj, v. 37, 1994,\\ \phantom{aaaaa} pp. 345 - 361.

9. Grigorian G. A. On two comparison tests for second-order linear  ordinary differential\\ \phantom{aaaaa} equations (Russian) Differ. Uravn. 47 (2011), no. 9, 1225 - 1240; translation in Differ.\\ \phantom{aaaaa} Equ. 47 (2011), no. 9 1237 - 1252, 34C10.

10. Demidovich B. P.  Lectures on the mathematical  stability theory, Moscow, ''Nauka'',\\ \phantom{aaaaa} 1967.

11. Pekarkova E. Asymptotic Properties of Second Order Differential Equations With \\ \phantom{aaaaa} p - Laplacian (Dissertation). Brno, Masarik University, Mathematics and Statistics, \\ \phantom{aaaaa}2009.

12. Momeni M.,  Kurakis I.,  Mosley - Fard M.,  Shukla P. K. A Van der Pol - Mathieu\\ \phantom{aaaaa} equation for the dynamics of dust grain charge in dusty plasmas. J. Phys, A: Math.\\ \phantom{aaaaa} Theor., v. 40, 2007,  F473 - F481.

 \end{document}